\documentclass{amsart}

\pdfoutput=1 

\usepackage[english]{babel}
\usepackage{graphicx}
\usepackage[hidelinks]{hyperref}
\usepackage{amssymb}
\usepackage{amsmath, mathtools}
\usepackage[foot]{amsaddr}
\usepackage{amsfonts}
\usepackage{bbm}
\usepackage{xcolor}
\usepackage{tikz}
\usepackage[top=1in, bottom=1.25in, left=1.25in, right=1.25in]{geometry}
\usepackage[bf, small]{caption}
\usepackage{subcaption}
\usepackage{multirow}
\usepackage{algorithm}
\usepackage{algpseudocode}
\usepackage{siunitx}

\definecolor{brightlavender}{rgb}{0.75, 0.58, 0.89}
\definecolor{DEblue}{rgb}{0.25, 0.0, 1.0}

\definecolor{maria}{HTML}{0090A0}

\definecolor{fabio}{HTML}{FF0000}

\theoremstyle{theorem} 
\newtheorem{Theorem}{Theorem}[subsection]
\newtheorem{Problem}[Theorem]{Problem}
\newtheorem*{Problem*}{Problem}
\newtheorem{Definition}[Theorem]{Definition}

\newtheorem{Assumption}[Theorem]{Assumption}

\newtheorem{Remark}[Theorem]{Remark}

\begin{document}
\title[SUPG ROMs for Advection-Dominated PDEs under Optimal Control]{A Streamline upwind Petrov-Galerkin Reduced Order Method for Advection-Dominated Partial Differential Equations under Optimal Control}
\author{Fabio Zoccolan$^1$, Maria Strazzullo$^2$, and Gianluigi Rozza$^3$}
\address{$^1$ Section de Mathématiques, École Polytechnique Fédérale de Lausanne, 1015 Lausanne, Switzerland, \newline email: fabio.zoccolan@epfl.ch}
\address{$^2$ DISMA, Politecnico di Torino, Corso Duca degli Abruzzi 24, Turin, Italy. \newline email: maria.strazzullo@polito.it}
\address{$^3$ mathLab, Mathematics Area, SISSA, via Bonomea 265, I-34136 Trieste, Italy. \newline email: gianluigi.rozza@sissa.it}

\begin{abstract}
In this paper we will consider distributed Linear-Quadratic Optimal Control Problems dealing with Advection-Diffusion PDEs for high values of the P\'eclet number. In this situation, computational instabilities occur, both for steady and unsteady cases. A Streamline Upwind Petrov–Galerkin technique is used in the optimality system to overcome these unpleasant effects. We will apply a finite element method discretization in an \textit{optimize-then-discretize} approach. {Concerning the parabolic case, a stabilized space-time framework will be considered and stabilization will also occur in both bilinear forms involving time derivatives}. Then we will build Reduced Order Models on this discretization procedure and two possible settings can be analyzed: whether or not stabilization is needed in the online phase, too. In order to build the reduced bases for state, control, and adjoint variables we will consider a Proper Orthogonal Decomposition algorithm in a partitioned approach.  {It is the first time that Reduced Order Models are applied to stabilized parabolic problems in this setting.}
The discussion is supported by computational experiments, where relative errors between the FEM and ROM solutions are studied together with the respective computational times.
\end{abstract}

\keywords{reduced order methods, time dependent parametrized optimal control problem, stabilization, proper orthogonal decomposition}

\maketitle

\section{Introduction}
\label{sec:intro}

The main goal of \textit{Optimal Control} theory is to modify a physical or engineering system through an input, called \textit{control}, to obtain a desired \textit{output}. From a theoretical point of view, one can describe the state problem through partial differential equations (PDEs), following the approach of J.L. Lions \cite{lions1971optimal,lions1972some}. Applying an optimal control means to solve a constrained optimization problem, where a cost functional has to be minimized. This process translates into an optimality system, which will be discretized for numerical simulations, that, in this framework, are more and more needed. Thus, effective and fast numerical techniques are required to exploit optimal control in scientific and industrial applications.

In this work, we will consider Advection-Diffusion equations \cite{quarteroni2008numerical} for large P\'eclet numbers. These equations are widespread in many engineering contexts since they can model transfer of particles, of energy, of heat and so on. In the case of high values of the P\'eclet number, numerical instabilities occur during discretization. We exploit a Streamline Upwind Petrov–Galerkin (SUPG) technique over a finite element method (FEM) \cite{brooks1982streamline, hughes1987recent, quarteroni2009numerical} in an \textit{optimize-then-discretize} approach, as done in \cite{collis2002analysis}, to provide strongly-consistency to the discretization. When we deal with parabolic problems, a space-time discretization  \cite{guberovic2014space,stoll2010all,strazzullo2020pod,strazzullo2021certified,strazzullo2022pod, urban2012new} will be used together with the SUPG stabilization for bilinear forms related to the derivative over time. {It is worth noticing that the use of a space-time discretization induces a solution of the system via a one-shot approach. The main advantage of this technique is to interpret the discretized parabolic problem as a steady one. In this way classical ROMs can be used with a small effort in generalization \cite{strazzullo2020pod}. However, as a drawback, increasing both time and space resolutions can lead to an unfeasible system to be solved in a many-query scenario. Therefore, reduction techniques can be a huge asset to alleviate this computational burden. We remark that other approaches can be used to solve the OCP system, such as iterative methods that have been used both at steady and unsteady levels. The interested reader may refer to \cite{manzoni2021optimal} for details. We stress that even these alternative approaches might only lighten the offline computational costs related to the problem formulation. Thus, for the sake of simplicity, we decided to work with the direct solution of the one-shot system as did previously in \cite{BALLARIN2022307,strazzullo2020pod,strazzullo2022pod}}.
In this case, the discretization procedure can easily request a huge amount of computational resources, especially for {parametric} time-dependent problems. The parameters can represent physical or geometrical features of the system at hand. 

In this scenario, we decide to exploit the parameter dependence of the equations to build Reduced Order Models (ROMs) \cite{hesthaven2016certified,quarteroni2014reduced,quarteroni2015reduced,rozza2008reduced} by means of Proper Orthogonal Decomposition (POD) algorithm in a partitioned approach. {Namely, the discretization process is divided in two phases: an \emph{offline} stage where a low-dimensional space is built through FEM solutions computed in properly chosen parameters, and an \emph{online} stage, where the system is solved for a new parametric instance in the new low-dimensional framework}. Thus, we consider two possible strategies: the former is to stabilize the system only in the offline phase; the latter uses SUPG in the online one, too.  {A priori, it is not always clear if instabilities can occur when the stabilization is settled only offline, as one could expect a "stabilizing effect" from the control itself, too. However, we will see that it is not the case in our numerical experiments and, hence, it becomes necessary to introduce online stabilization to overcome this undesired behaviour on the reduced solution.} This setting was considered for problems without optimal control in \cite{pacciarini2014stabilized,torlo2018stabilized}. {Moreover, to the best of our knowledge, it is the first time that SUPG stabilization for time-dependent Advection-Dominated problems under distributed control is analyzed in a ROM setting}.

This work is organized as follows.
The first section will illustrate some theoretical aspects about Optimal Control Theory for PDEs. Section \ref{sec:truth_discretization} shows an introduction to Advection-Dominated problems and the SUPG technique over FEM discretization in an \textit{optimize-then-discretize} approach. In Section \ref{sec:ROMs}, we will focus on the ROM setting and Section \ref{sec:results} refers to the related numerical simulations. Firstly, we will introduce two specific examples of Advection-Diffusion problems: the Graetz-Poiseuille and the Propagating Front in a Square Problems. The former was studied in various forms without optimal control in \cite{gelsomino2011comparison,pacciarini2014stabilized,rozza2009reduced,torlo2018stabilized} and with optimal control but without stabilization in \cite{negri2011reduced,strazzullo2020pod}. The latter is studied without optimal control in a similar version in \cite{pacciarini2014stabilized,torlo2018stabilized}. Here, both the problems will be analyzed under a distributed optimal control for high values of the P\'eclet number, both in the steady and unsteady cases. Relative errors between FEM and ROM solutions will be shown, as well as an analysis on the computational times.

\section{Problem Formulation}
\label{sec:problem}

{In this Section we will illustrate the basic properties of Linear-Quadratic Optimal Control Problem (OCP) for steady and parabolic PDEs. The interested reader may refer to \cite{doi:10.1137/090760453,gelsomino2011comparison,karcher2018certified,negri2015reduced,negri2013reduced,quarteroni2008numerical} for steady problems and to \cite{strazzullo2018model,strazzullo2020pod,strazzullo2021certified} for parabolic ones for a detailed discussion. In this work we will only consider {\textit{distributed control}, which means that the control acts on the whole physical domain.}}

{The aim of Optimal Control is to achieve a {prescribed} optimality condition by minimizing a suitable cost functional under the constraint of satisfying  the PDE Problem. {The proposed framework} follows the J.L. Lions theory \cite{lions1971optimal, lions1972some}. From a mathematical perspective, we can state that an OCP is characterized by:
\begin{itemize}
    \item $e$, the \textit{state equation function}, which expresses the relationship between the output and the control within the system in terms of a PDE problem or PDEs in a weak formulation. A pair $(y,u) \in X:=Y\times U$ is said to be $physical$ or $feasible$ if it is a solution of the state equation $e$; $y$ is called the \textit{state} variable, \textit{the output}, and $u$ is the \text{control} variable, \textit{the input}. $X_{ad}$ is the set of all the feasible pairs $(y,u)$;
    \item $z(y) = Oy$, a direct \textit{observation} of the output. Here, a linear operator $O$ is applied to the state to describe the observation: we will denote the space of observation as $Z$. We will only deal with state variables that can be measured on a portion of the domain;
    \item $J$, the \textit{objective functional}, which describes the objective to achieve;
    \item {suitable spaces $Y$ and $U$, as the \textit{state space} and \textit{control space}, respectively. Domains of definition for control and/or state can be taken smaller due to possible restrictions; hence we introduce $Y_{ad} \subseteq Y$ and $U_{ad} \subseteq U$ as the \textit{admissible state space} and \textit{control space}, respectively. However, we will always consider \textit{unconstrained problems}, i.e. $X_{ad} = X$. The theory of the well-posedness can make use of the Lagrangian approach \cite{carere2019thesis, negri2011reduced} or of the general Adjoint one when dealing with $X_{ad} \subset Y_{ad} \times U_{ad}$ \cite{hinze2008optimization, lions1971optimal, quarteroni2009numerical}.}
\end{itemize}
}
Let us consider $\Omega \subset \mathbb{R}^{n}$, an open and bounded regular domain, and the time interval $(0, T) \subset \mathbb{R}^{+}$: for us, it will always be the case of $n=2$. Let us denote with $\Gamma_{D}$ and $\Gamma_{N}$ the portions of the boundary of $\partial \Omega$ where Dirichlet and Neumann boundary conditions are specified, respectively. {We define the \textit{observation domain} $\Omega_{obs} \subseteq \Omega$ as the portion of the domain where we want the state variable to assume a desired value.} $\mathcal{P} \subseteq \mathbb{R}^p$, for natural number $p$, is the parameter space and $\mathbf{\boldsymbol{\mu}} \in \mathcal{P}$ is a $p$-vector which can represent physical or geometrical parameter of interest.  In this work, we deal with \textit{Parametric Optimal Control Problems} (OCP($\boldsymbol \mu$)s), i.e.\ systems where there is a dependency on the parameter $\boldsymbol{\mu}$.

\subsection{Steady Problems}
\label{sec-Lagrangian}
We refer to the Lagrangian approach to state the well-posedness of OCP($\boldsymbol \mu$)s in \textit{full admissibility} setting, i.e.\ when $X_{ad} = Y \times U$, \cite{carere2019thesis,hinze2008optimization,troltzsch2010optimal}. We want to solve:
\begin{equation}\label{OCP-general}
    \min_{(y(\boldsymbol{\mu}),u(\boldsymbol{\mu})) \in Y \times U} J(y(\boldsymbol{\mu}),u(\boldsymbol{\mu}); \boldsymbol{\mu}) \ \text{s.t.} \ e(y(\boldsymbol{\mu}),u(\boldsymbol{\mu}); \boldsymbol{\mu})= 0,
\end{equation}
thus we define the Lagrangian operator $\mathcal{L} : Y \times U \times Q^{*} \to \mathbb{R}$ as:
\begin{equation}
  \label{lagrangian}
   \mathcal{L}(y(\boldsymbol{\mu}),u(\boldsymbol{\mu}),p(\boldsymbol{\mu}); \boldsymbol{\mu}) = J(y(\boldsymbol{\mu}),u(\boldsymbol{\mu}); \boldsymbol{\mu})+\langle p(\boldsymbol{\mu}),e(y(\boldsymbol{\mu}),u(\boldsymbol{\mu}); \boldsymbol{\mu})\rangle_{Q^{*}Q},
\end{equation}
where {$p(\boldsymbol{\mu})$ is a Lagrange multiplier belonging to} $Q^{*}$, the dual space of $Q$. For the sake of notation, we write $y:=y(\boldsymbol{\mu})$, $u:=u(\boldsymbol{\mu})$ and $p:=p(\boldsymbol{\mu})$: we will explicit the parameter dependence only when necessary. 
{We remark that in this paper we will always consider Linear-Quadratic problems.}
\begin{Problem}[Linear-Quadratic Problem]\label{def:lin-quad-prob}
{Consider the Banach spaces  $Y,U$, and $Z$, and $\alpha>0$. Let $\mathbf{\boldsymbol{\mu}} \in \mathcal{P}$ be given.} Let the Observation map $O : Y \to Z$ be a linear and bounded operator. Consider an element $z_{d}(\boldsymbol{\mu}) \in Z$, which is the so-called desired solution profile (the desired observed output). Let $J$ be a quadratic objective functional of the form
\begin{equation}
\label{quadratic-functional}
   J(y, u; \boldsymbol{\mu})=\frac{1}{2} m\left(O y(\boldsymbol{\mu})-z_{d}(\boldsymbol{\mu}), Oy(\boldsymbol{\mu})-z_{d}(\boldsymbol{\mu})\right)+\frac{\alpha}{2} n(u(\boldsymbol{\mu}), u(\boldsymbol{\mu})),
\end{equation}
where $m: Z \times Z \rightarrow \mathbb{R} $ and $n: U \times U \rightarrow \mathbb{R}$ are symmetric and continuous bilinear forms.
{Let $e : Y \times U \to Q$, with $Q$ a Banach space, which fulfills a set of boundary and/or initial conditions, }be affine, i.e.\ there exist $A(\boldsymbol{\mu}) \in \mathcal{B}(Y, Q)$, $B(\boldsymbol{\mu}) \in \mathcal{B}(U, Q)$ and $f(\boldsymbol{\mu}) \in Q$ such that
\begin{equation}
\label{affine-state}
e(y, u; \boldsymbol{\mu})=A(\boldsymbol{\mu}) y+B(\boldsymbol{\mu}) u-f(\boldsymbol{\mu}), \quad \forall \big(y(\boldsymbol{\mu}), u(\boldsymbol{\mu})\big) \in Y \times U.
\end{equation}
Then {the minimization problem} \eqref{OCP-general} with the above properties is said a \textit{Linear-Quadratic Optimal Control Problem}.
\end{Problem}
{Under the assumption that $A(\boldsymbol{\mu})$ admits bounded inverse, Problem \ref{def:lin-quad-prob} has a unique solution \cite[Theorem 1.43]{hinze2008optimization}. Moreover,
For Linear-Quadratic OCP($\boldsymbol \mu$)s,}{a pair $(y, u)$ to Problem \ref{OCP-general} is a solution if and only if there exists a $p \in Q^{*}$ such that \eqref{Opt-lin_system} holds \cite{carere2019thesis,troltzsch2010optimal}: }
\begin{equation} \label{Opt-lin_system}
\begin{aligned}
\begin{cases}
m(O y, O \bar{y}; \boldsymbol{\mu})+\left\langle A^{*}(\boldsymbol{\mu}) p, \bar{y}\right\rangle_{Y^{*} Y} =m\left(O \bar{y}, z_{d}; \boldsymbol{\mu}\right), & \forall \bar{y} \in Y, \\
\alpha n(u, \bar{u}; \boldsymbol{\mu})+\left\langle B^{*}(\boldsymbol{\mu}) p, \bar{u}\right\rangle_{U^{*} U} =0, & \forall \bar{u} \in U, \\
\langle\bar{p}, A(\boldsymbol{\mu}) y+B(\boldsymbol{\mu}) u\rangle_{Q^{*} Q} =\langle\bar{p}, f(\boldsymbol{\mu})\rangle_{Q^{*} Q}, & \forall \bar{p} \in Q^{*}.
\end{cases}
\end{aligned}
\end{equation}
In \eqref{Opt-lin_system}, the first equation is called the \textit{adjoint equation}, the second one is the \textit{gradient equation}, whereas the \textit{state equation} is the third one.

\begin{Remark}[Notation] \label{remark-notation}
  From now on, we always involve Hilbert spaces. For the sake of notation, we will denote the various bilinear forms defined by $A(\boldsymbol{\mu}),B(\boldsymbol{\mu})$ in the following unique way:
  
  \begin{equation*}
      \begin{aligned}
           \langle A(\boldsymbol{\mu}) y, p\rangle_{Q Q^{*}} := a(y,p;\boldsymbol{\mu}) \qquad & \langle B(\boldsymbol{\mu}) u, p\rangle_{Q Q^{*}} := b(u,p; \boldsymbol{\mu}).
      \end{aligned}
  \end{equation*}
\end{Remark}

\subsection{Unsteady Problems}
\label{Unsteadyocp}
We briefly introduce time-dependent Linear-Quadratic OCP($\boldsymbol \mu$)s based on \cite{strazzullo2020pod,strazzullo2021certified}. One can consider saddle-point formulation in order to prove well-posedness by using tools of the previous Sections in the case of null initial conditions.
Differently from the steady case, here one needs some more technical assumptions, which will be fulfilled by both “Graetz-Poiseuille" and “Propagating Front in a Square" problems.

Consider two separable Hilbert spaces $Y$ and $H$ satisfying $Y \hookrightarrow H \hookrightarrow Y^{*}$ and, moreover, other two Hilbert spaces $U$ and $Z \supseteq Y$, where $Y$ and $U$ are the usual \textit{state} and \textit{control spaces}, and $Z$ is the space of observation. We endow them with the standard norms inherited from their respectively scalar products: $(\cdot,\cdot)_{Y}$, $(\cdot,\cdot)_{Z}$, $(\cdot,\cdot)_{U}$, and $(\cdot,\cdot)_{H}$.
We define the following Hilbert spaces $\mathcal{Y}=L^{2}(0, T ; Y),$ $\mathcal{Y}^{*}= L^{2}\left(0, T ; Y^{*}\right)$, $\mathcal{U}=L^{2}(0, T ; U)$, $\mathcal{Z}:=L^{2}(0, T ; Z) \supseteq \mathcal{Y},$
with respective norms, for instance in the case of $\mathcal{Y}$ and $\mathcal{U}$ given by
$\|y\|_{\mathcal{Y}}^{2}:=\int\limits_{0}^{T}\|y\|_{Y}^{2} \mathrm{dt}$, and  $\|u\|_{\mathcal{U}}^{2}:=\int\limits_{0}^{T}\|u\|_{U}^{2} \mathrm{dt}$,
and similarly for the others.
Furthermore, let us define the Hilbert space $\mathcal{Y}_{t}$  with its scalar product $(\cdot, \cdot)_{\mathcal{Y}_{t}}$:
\begin{equation*} \label{y_t}
    \mathcal{Y}_{t}:=\left\{ y \in \mathcal{Y} \quad \text{s.t.} \quad \frac{\partial y}{\partial t} \in \mathcal{Y}^{*} \right\}, \quad
\displaystyle
(y,z)_{\mathcal{Y}_{t}}:=\int\limits_{0}^{T}(y,z)_{Y} \mathrm{dt}+\int\limits_{0}^{T}\Big(\frac{\partial y}{\partial t}, \frac{\partial z}{\partial t}\Big)_{Y^{*} } \mathrm{dt}.
\end{equation*}
Our aim is to solve the following unconstrained Linear-Quadratic \textit{Parametric Parabolic OCP$(\boldsymbol{\mu})$}:
\begin{Problem}[Parametric Parabolic OCP$(\boldsymbol{\mu})$]\label{Parametrized-Parabolic-OCP}
   Given $\boldsymbol{\mu} \in \mathcal{P}$ find the pair $(y(\boldsymbol{\mu}), u(\boldsymbol{\mu})) \in \mathcal{Y}_{t} \times \mathcal{U}$ 
satisfying
\begin{equation}\label{par-prob}
\begin{cases}
\displaystyle
\frac{\partial y(\boldsymbol{\mu})}{\partial t}+A(\boldsymbol{\mu}) y(\boldsymbol{\mu})+B(\boldsymbol{\mu}) u(\boldsymbol{\mu})-f(\boldsymbol{\mu})=0, & \text { in } \Omega \times(0, T), \\
\displaystyle
\frac{\partial y(\boldsymbol{\mu})}{\partial n}=0, & \text { on } \Gamma_{N} \times(0, T), \\
\displaystyle
y(\boldsymbol{\mu})=l, & \text { on } \Gamma_{D} \times(0, T), \\
\displaystyle
y(\boldsymbol{\mu})(0)=y_{0}, & \text { in } \Omega,
\end{cases}
\end{equation}
$$
\mbox{and } \qquad \min _{(y(\boldsymbol{\mu}), u(\boldsymbol{\mu})) \in \mathcal{Y}_{t} \times \mathcal{U}} J(y, u ; \boldsymbol{\mu}) = \frac{1}{2} m\left(O y(\boldsymbol{\mu})-z_{d}(\boldsymbol{\mu}), O y(\boldsymbol{\mu})-z_{d}(\boldsymbol{\mu})\right)+\frac{\alpha}{2} n(u(\boldsymbol{\mu}), u(\boldsymbol{\mu})), 
$$
where $m: \mathcal{Y}_{t} \times \mathcal{Y}_{t} \rightarrow \mathbb{R} $ and $n: \mathcal{U} \times \mathcal{U} \rightarrow \mathbb{R}$ are symmetric and continuous bilinear forms, $z_{d}({\boldsymbol{\mu}}) \in \mathcal{Z}$ is the observed desired solution profile and $\alpha>0$ is the fixed penalization parameter. In our test case we will always take $y_0 \equiv 0$.
\end{Problem}

Again, we denote $y:=y(\boldsymbol{\mu})$ and $u:=u(\boldsymbol{\mu})$ omitting the parameter dependence. 
{Assuming that $A(\boldsymbol{\mu})$ admits a bounded inverse and $f(\boldsymbol{\mu}) \in \mathcal{Y}^{*}$, then there exists a unique solution for Problem \ref{Parametrized-Parabolic-OCP} \cite{hinze2008optimization,troltzsch2010optimal}. Moreover, if we assume that $f(\boldsymbol{\mu})$ }gathers all forcing, boundary {and, eventually, lifting} terms of the state equation, we can state a Lagrangian approach based on the weak formulation of Problem \ref{Parametrized-Parabolic-OCP}. 
For the sake of notation, we denote $a: \mathcal{Y}_{t} \times \mathcal{Y}_{t} \rightarrow \mathbb{R}$ and $b: \mathcal{U} \times \mathcal{Y}_{t} \rightarrow \mathbb{R}$ the bilinear forms defined as $a(y,q; \boldsymbol{\mu})= \langle A(\boldsymbol{\mu})y, q \rangle_{\mathcal{Y}^{*}\mathcal{Y}}$ and $b(u,q; \boldsymbol{\mu})=  \langle B(\boldsymbol{\mu})u, q \rangle_{\mathcal{Y}^{*}\mathcal{Y}}$, respectively. For a proper definition of the adjoint variable, it is opportune to take $q \in \mathcal{Y}_{t}$ rather than $q \in \mathcal{Y}$ \cite{strazzullo2020pod}. Denoting $p:=p(\boldsymbol{\mu})$ the adjoint test function and considering $Q^{*}=\mathcal{Y}_{t}$,  in our setting $p(T)=0$ holds \cite{lions1971optimal}.
\begin{Remark} \label{null-initial and Z}
   For both steady and unsteady problems, we take{the embedding between $Y=H^{1}(\Omega)$ and $Z=L^2(\Omega_{obs})$} as the Observation function $O$. Thus, 
   $\mathcal{Z}=\mathcal{Y}$ and our desired state is denoted by $y_d$.
\end{Remark}
{In the next section we extend the proposed formulation to the Advection-Dominated setting.}\color{black}
\section{Truth Discretization}
\label{sec:truth_discretization}

{In this {Section} we present the involved discretization procedure. First, we will consider SUPG stabilization for Advection-Dominated equations in case of high Péclet number. Then, an \textit{optimize-then-discretize} approach is followed, i.e.\ at first we derive optimality conditions for the stabilized PDE as system \eqref{Opt-lin_system} and then we discretize it following an \textit{one shot} or \textit{all-at-once approach} \cite{hinze2008hierarchical,stoll2010all, STOLL2013498}. }

Let us start our discussion from the steady case. From now on we will always assume to work with $Y,U,Q$ Hilbert spaces. We employ a FEM discretization, which will be named as the \textit{high-fidelity} or \textit{truth} approximation. We consider $\Omega_h$ as a quasi-uniform mesh on the domain $\Omega$, where $h$ indicates the \textit{mesh size}, i.e. maximum diameter of an element of the grid. $\mathcal{T}_h$ is {a regular triangularization} on $\Omega$ and
$$
\Omega_{h}:=\operatorname{int}\left(\bigcup_{K \in \mathcal{T}_{h}}{K}\right),
$$
where $K$ is a triangle of $\mathcal{T}_{h}$. {We define the FEM spaces $Y^{\mathcal{N}}=Y \cap X^{\mathcal{N}, 1}$, $U^{\mathcal{N}}=U \cap X^{\mathcal{N}, 1}$ and $\big(Q^{\mathcal{N}}\big)^{*}={Q^{*}} \cap X^{\mathcal{N}, 1}$, where
$$
X^{\mathcal{N}, 1}=\left\{v^{\mathcal{N}} \in C^{0}(\bar{\Omega}): v^{\mathcal{N}}_{|_{K}} \in \mathbb{P}^{1}(K), \forall K \in \mathcal{T}_{h}\right\}
$$
and $\mathbb{P}^{1}(K)$ represents the space of polynomials of degree at most equal to $1$ defined on $K$. We will always use the same triangulation $\mathcal{T}_{h}$ and a $\mathbb{P}^{1}$-FEM approximation for state, control and adjoint variables. The dimensions of $Y^{\mathcal{N}}, U^{\mathcal{N}}, Q^{\mathcal{N}}$ are all equal to $\mathcal{N}$. The overall dimension of the discrete problem is $\mathcal{N}_{tot} = 3 \cdot \mathcal{N}$. For the sake of simplicity, we assume ${Q}^{*}_{h} = Y^{\mathcal{N}}$. 
Let us denote with $\mathbf{y}$, $\mathbf{u}$ and $\mathbf{p}$ the {vectors of coefficients} of $y^{\mathcal{N}}$, $u^{\mathcal{N}}$ and $p^{\mathcal{N}}$, expressed in terms of the nodal basis of $Y^{\mathcal{N}}$, $U^{\mathcal{N}}$, and $(Q^{\mathcal{N}})^{*}$. We express with $\boldsymbol{y}_d$ the vector with the components of the discretized desired state, i.e.\ the Galerkin projection of $y_d$ on $Y^{\mathcal{N}}$. Moreover, from now on the operators $m$ and $n$ in \eqref{quadratic-functional} will be the $L^2$ product on $\Omega_{obs}$ and on $\Omega$, respectively. }

\subsection{SUPG stabilization for Advection-Dominated OCP($\boldsymbol{\mu}$)s}
In this Section we illustrate Advection-Dominated OCP($\boldsymbol{\mu}$)s and the SUPG technique applied to an \textit{optimize-then-discretize} approach. From now, we recall the dependence on parameters of our operators. Let us start from our definition of an \textit{Advection-Diffusion equation}.

\begin{Definition}[Advection-Diffusion Equations]\label{advection-diffusion}
 Let us consider the following problem:
\begin{equation}\label{advection-diffusion problem}
L (\boldsymbol{\mu})y:=-\varepsilon(\boldsymbol{\mu}) \Delta y+\mathbf{b} (\boldsymbol{\mu}) \cdot \nabla y=f(\boldsymbol{\mu})  \text { in } \Omega \subset \mathbb{R}^2,
\end{equation}
with suitable boundary conditions on $ \partial \Omega$. Let us suppose that:
 \begin{itemize}
     \item the diffusion coefficient $\varepsilon: \Omega \rightarrow \mathbb{R}$ belongs to $L^{\infty}(\Omega)$ and depends on the parameter $\boldsymbol{\mu}$. We assume there exists a constant $\bar{\varepsilon}>0$ such that
$
\varepsilon(x) \geq \bar{\varepsilon}$, for almost all $x \in \Omega;
$
\item the advection field $\boldsymbol{b}: \Omega \rightarrow \mathbb{R}^{2}$ belongs to $\left(L^{\infty}(\Omega)\right)^{2}$ and depends on the parameter $\boldsymbol{\mu}$. {We suppose that 
$\operatorname{div} \mathbf{b}(x)$ exists and $\operatorname{div} \mathbf{b}(x)=0$ holds for almost all $x \in \Omega$};
\item $f(\boldsymbol{\mu}): \Omega \rightarrow \mathbb{R}$ is an $L^{2}(\Omega)$-function that can depend on the parameter $\boldsymbol{\mu}$.
 \end{itemize}
 
In this case, \eqref{advection-diffusion problem} is an Advection-Diffusion problem and $L (\boldsymbol{\mu})$ is the Advection-Diffusion operator.
\end{Definition}
{From \eqref{advection-diffusion problem}, we can easily derive the weak formulation of an Advection-Diffusion problem: find $y \in  Y$ s.t. 
$a\left(y, q ; \boldsymbol{\mu}\right)=F\left(q ; \boldsymbol{\mu}\right) \text{ for all } q \in Q^{*}$, where 
\begin{equation}\displaystyle
     a\left(y, q ; \boldsymbol{\mu}\right) := \int_{\Omega} \varepsilon(\boldsymbol{\mu}) \nabla y \nabla q + \boldsymbol{b}(\boldsymbol{\mu}) \cdot \nabla y q  \ \mathrm{dx}, \quad F(q; \boldsymbol{\mu}) :=\int_{\Omega} f(\boldsymbol{\mu}) q \ \mathrm{dx}, \quad y \in Y,q \in Q^{*}.
\end{equation}


From a numerical point of view, when the advection term $\mathbf{b}(\boldsymbol{\mu})\cdot \nabla u$  ``dominates'' the diffusive one $-\varepsilon (\boldsymbol{\mu}) \Delta u$, i.e. when $|\mathbf{b}(\boldsymbol{\mu})| \gg \varepsilon(\boldsymbol{\mu})$, the approximated solution can show instability phenomena along the direction of the advection field \cite{quarteroni2008numerical}.
In order to give an indicator of the instability, let us consider the regular triangulation $\mathcal{T}_{h}$ related to FEM discretization. For any element $K \in \mathcal{T}_{h}$, where $h_{K}$ is the diameter of $\mathrm{K}$, we can then define the \textit{local P\'eclet number} as \cite{quarteroni2008numerical, quarteroni2009numerical}:
\begin{equation}\label{Péclet number}
\mathbb{P}\mathrm{e}_{K}(x):=\frac{|\mathbf{b}(x)| h_{K}}{2 \varepsilon(x)}, \quad \forall x \in K.
\end{equation}

\begin{Definition}[Advection-Dominated problem]\label{adv-problem-def}
 Considering Definition \ref{advection-diffusion}, we are dealing with an Advection-Dominated problem if
 $
\mathbb{P}\mathrm{e}_{K}(x)>1, \ \forall x \in K, \ \forall K \in \mathcal{T}_{h}$ {\cite{quarteroni2008numerical}}.
\end{Definition}

To solve the issue of the instability, we will exploit the SUPG method \cite{brooks1982streamline,hughes1979multidimentional,hughes1987recent, quarteroni2008numerical}, which is a \textit{strongly consistent stabilization technique}; i.e.\ it is consistent for weak PDEs and its order of accuracy can be greater than one.
Let us now consider the Advection-Diffusion operator \eqref{advection-diffusion problem}: for the sake of simplicity, we define it on $H_{0}^{1}(\Omega)$ and we do not indicate the parameter dependence.
The operator $L$ can be split into its symmetric and skew-symmetric parts \cite{quarteroni2008numerical}, defined as:
\begin{equation}\label{sym-skewsym}
    \text{symmetric part: }  L_{S} y=-\varepsilon \Delta y, \qquad
    \text{skew-symmetric part: }  L_{S S} y=\mathbf{b} \cdot \nabla y; \quad \forall y \in H_{0}^{1}(\Omega),
\end{equation}
i.e.\ $L=L_{S}+L_{S S}.$ Symmetric and skew-symmetric parts can be directly derived using the formulae:
\begin{equation}\label{sks-s-parts}
L_{S} =\frac{L+L^{*}}{2}, \qquad L_{S S} =\frac{L-L^{*}}{2},
\end{equation}
where $L^{*}$ is the adjoint operator related to $L$. 
Now, let us analyze our OCP \eqref{Opt-lin_system}: we follow the \textit{optimize-then-discretize} approach in \cite{collis2002analysis}. The \textit{discretized state equation} is described as follows, where the control is distributed, i.e.\ it acts on the whole domain $\Omega$: 
\begin{equation}\label{discretized-se}
    a_{s}\left(y^{\mathcal{N}}, q^{\mathcal{N}}\right)+b_{s}\left(u^{\mathcal{N}},q^{\mathcal{N}}\right)=F_{s}(q^{\mathcal{N}}), \quad \forall q^{\mathcal{N}} \in \left(Q^{\mathcal{N}}\right)^{*},
\end{equation}
with
\begin{equation}\label{supg-form}
     a_{s}\left(y^{\mathcal{N}}, q^{\mathcal{N}}\right):=a\left(y^{\mathcal{N}}, q^{\mathcal{N}} \right)+\sum_{K \in \mathcal{T}_{h}} \delta_{K}\left(L y^{\mathcal{N}}, \frac{h_{K}}{|\mathbf{b}|} L_{SS}q^{\mathcal{N}}\right)_{K}, 
\end{equation}
\begin{equation}
    b_{s}\left(u^{\mathcal{N}},q^{\mathcal{N}}\right):= - \int\limits_{\Omega} u^{\mathcal{N}} q^{\mathcal{N}} - \sum_{K \in \mathcal{T}_{h}} \delta_{K}\left( u^{\mathcal{N}}, \frac{h_{K}}{|\mathbf{b}|} L_{SS}q^{\mathcal{N}}\right)_{K},
\end{equation}
and 
\begin{equation}
F_{s}(q^{\mathcal{N}}) := F\left(q^{\mathcal{N}}\right) + \sum_{K \in \mathcal{T}_{h}} \delta_{K}\left( f, \frac{h_{K}}{|\mathbf{b}|} L_{SS}q^{\mathcal{N}}\right)_{K},
\end{equation}
where $\big(\cdot, \cdot \big)_K$ indicates the usual $L^2(K)$-product, $f$ collects all the forcing and lifting terms, and $\delta_{K}$ denotes a positive dimensionless stabilization parameter related to an element $K \in \mathcal{T}_{h}$. In principle, since $\delta_{K}$ is local, it can be different for each $K$. The adjoint equation is an Advection-Dominated equation, too, but with an advective term with opposite sign with respect to the state one. As a matter of fact, from \eqref{sks-s-parts} we obtain that $L^{*} = L_{S}-L_{SS}$. The SUPG method leads to the \textit{discretized adjoint equation}
\begin{equation}\label{discretized-ae}
a_{s}^{*}\left(z^{\mathcal{N}}, p^{\mathcal{N}}\right)+\big( y^{\mathcal{N}}-y_d, z^{\mathcal{N}}\big)_{s}=0, \quad \forall z^{\mathcal{N}} \in Y^{\mathcal{N}},
\end{equation}
with
\begin{equation}
\begin{aligned}
a_{s}^{*}\left(z^{\mathcal{N}}, p^{\mathcal{N}}\right) &:=a^{*}\left(z^{\mathcal{N}}, p^{\mathcal{N}}\right)+\sum_{K \in \mathcal{T}_{h}} \delta^{a}_{K}\left( (L_{S}-L_{SS}) p^{\mathcal{N}}, \frac{h_{K}}{|\mathbf{b}|}\left(-L_{SS}\right) z^{\mathcal{N}}\right)_{K}, \\
\big( y^{\mathcal{N}}-y_d, z^{\mathcal{N}}\big)_{s} &:= \int\limits_{\Omega_{obs}} (y^{\mathcal{N}}-
    y_d)z^{\mathcal{N}} \ \mathrm{dx} +\sum_{K \in {\mathcal{T}_{h}}_{\vert_{\Omega_{obs}}}} \delta^{a}_{K}\left( y^{\mathcal{N}}-y_d, \frac{h_{K}}{|\mathbf{b}|}\left(-L_{SS}\right) z^{\mathcal{N}}\right)_{K},
\end{aligned}
\end{equation}
where $a^{*}$ is the adjoint form of $a$ and $\delta^{a}_{K}$ is the parameter related to the stabilized adjoint bilinear forms. As in this work we consider $\delta_{K}=\delta^{a}_{K}$ in numerical simulations, from now on we will always denote both stabilization parameter by $\delta_{K}$. Instead, the \textit{discretized gradient equation} remains untouched:
\begin{equation}
b^{*}\big(v^{\mathcal{N}}, p^{\mathcal{N}}\big)+\alpha n\big(u^{\mathcal{N}}, v^{\mathcal{N} }\big)=0, \quad \forall v^{\mathcal{N}} \in U^{\mathcal{N}}.
\end{equation}
With this setting it follows a nonsymmetric system for the computation of the numerical solution, but we gain the strong consistency of the method for the optimality system if $y, u, p$ are regular \cite{collis2002analysis}.
To summarize, the SUPG optimality system for a steady OCP is the following:
\begin{equation}\label{supg-system}
    \begin{aligned}
    &\textit{discretized adjoint equation:} \quad &a_{s}^{*}\left(z^{\mathcal{N}}, p^{\mathcal{N}}\right)+\big( y^{\mathcal{N}}-y_d, z^{\mathcal{N}}\big)_{s}=0 , \quad &\forall z^{\mathcal{N}} \in Y^{\mathcal{N}},\\
   &\textit{discretized gradient equation:} \quad &b^{*}\big(v^{\mathcal{N}}, p^{\mathcal{N}}\big)+\alpha n\big(u^{\mathcal{N}}, v^{\mathcal{N}}\big)=0 , \quad & \forall v^{\mathcal{N}} \in U^{\mathcal{N}}, \\
         &\textit{discretized state equation:} \quad &a_{s}\left(y^{\mathcal{N}}, q^{\mathcal{N}}\right)+b_{s}\left(u^{\mathcal{N}},q^{\mathcal{N}}\right)=F_{s}(q^{\mathcal{N}}) , \quad & \forall q^{\mathcal{N}} \in \left(Q^{\mathcal{N}}\right)^{*},
    \end{aligned}
\end{equation}
and the discretized algebraic system reads as:
\begin{equation} \label{stab-block-system}
\left(\begin{array}{ccc}
M_{s} & 0 & K_{s}^{T} \\
0 & \alpha M &B^{T} \\
K_{s} &B_{s} & 0
\end{array}\right)\left(\begin{array}{l}
\mathbf{y} \\
\mathbf{u} \\
\mathbf{p}
\end{array}\right)=\left(\begin{array}{c}
M_{s} \mathbf{y}_{d} \\
0 \\
\mathbf{f}_{s}
\end{array}\right),
\end{equation}
where $M_{s}$ is the stabilized mass matrix related to $m$, $M$ is the not-stabilized mass matrix related to $n$, $K_{s}$ and $K_{s}^{T}$ are the stiffness matrices related to $a_{s}$ and $a_{s}^{*}$, respectively, $B_s$ is the stabilized mass matrix related to $b_s$, $B^{T}$ is the block linked to $b^{*}$ and $\mathbf{f}_{s}$ is the vector whose components are the coefficients of the stabilized force term. 

We indicate with $|\|\cdot\||$ the energy norm of the bilinear form $a$ in the Advection-Diffusion operator, i.e. 
\begin{equation} \label{3I norm}
|\| w\||^{2}:=\varepsilon\|\nabla w\|_{L^{2}(\Omega)}^{2}+\frac{1}{2}\left\|(\operatorname{div} \boldsymbol{b})^{\frac{1}{2}} w\right\|_{L^{2}(\Omega)}^{2}, \quad \forall w \in H_{0}^{1}(\Omega).
\end{equation}
Therefore, we define the SUPG norm on $H_{0}^{1}(\Omega)$ as
\begin{equation} \label{supg-norm}
\|w\|_{S U P G}^{2}:=|\|w\||^{2}+\sum_{K \in \mathcal{T}_{h}} \delta_{K}\left(L_{S S} w, \frac{h_{K}}{|\boldsymbol{\boldsymbol{b}}|} L_{S S} w\right)_{K}, \quad \forall w \in H_{0}^{1}(\Omega).
\end{equation}
{Under local assumptions on the stabilization parameter $\delta_K$ (we postpone their definition in Theorem \ref{error-ocp-supg}), the SUPG bilinear form \eqref{supg-form} can be proved to be coercive with respect to the SUPG norm \cite{quarteroni2008numerical}, which implies the well-posedness of the stabilized problem.}
Finally, we can illustrate an estimate of the error for the adjoint and the state variables of the solution of an OCP \cite{collis2002analysis}{: we would like to point out that in \cite{heinkenschloss2010local} a refined error local analysis in this context is also presented.}
\begin{Theorem}[Error for state and adjoint variables]\label{error-ocp-supg}
    Let $m, r \geq 1$ and $(y, u, p)$ be the solution of \eqref{Opt-lin_system} with $y \in H^{m+1}(\Omega),  p \in H^{r+1}(\Omega) .$ Furthermore, let $y^{\mathcal{N}}, u^{\mathcal{N}}, p^{\mathcal{N}}$ be the numerical solution of \eqref{supg-system}.  If $\delta_K$ satisfies
    \begin{equation} \label{delta-req}
        0 < \delta_K \leq \frac{h_K}{\varepsilon \eta^{2}_{\text{inv}}} \quad \text{and} \quad \delta_K = \begin{cases}
                     \begin{aligned}
           \delta_1 \frac{h_{K}}{\varepsilon}, \quad  & \mathbb{P}\mathrm{e}_{K}(x) \leq 1, \\
           \delta_2, \quad & \mathbb{P}\mathrm{e}_{K}(x) > 1, 
                      \end{aligned}
             \end{cases}
    \end{equation}
    where $\delta_1, \delta_2 >0$ are chosen constant, 
    and $\eta_{\text{inv}}$ is defined as the following inverse constant
    \begin{equation*}
        |y^{\mathcal{N}}|_{1,K} \leq \eta_{\text{inv}}h^{-1}_{K}\|y^{\mathcal{N}}\|_{L^2(K)} \quad \text{and} \quad  \|\Delta y^{\mathcal{N}}\|_{L^2(K)} \leq \eta_{\text{inv}}h^{-1}_{K}\|\nabla y^{\mathcal{N}}\|_{L^2(K)} \quad \forall y^{\mathcal{N}} \in Y^{\mathcal{N}},
    \end{equation*} with $|\cdot|_{1,K}$, $\|\cdot\|_{K}$ seminorm and $L^2$-norm on $K$, respectively, then there exists $C>0$ such that
\begin{equation}
\begin{aligned}
\left\|y-y^{\mathcal{N}}\right\|_{SUPG} & \leq C\left(h^{m}\left(\varepsilon^{1 / 2}+h^{1 / 2}\right)|y|_{H^{m+1}(\Omega)} +\left\|u^{\mathcal{N}}-u\right\|_{L^2(\Omega)}\right), & \forall h, \\
\left\|p-p^{\mathcal{N}}\right\|_{SUPG} & \leq C\left(h^{r}\left(\varepsilon^{1 / 2}+h^{1 / 2}\right)|p|_{H^{r+1}(\Omega)} +\left\|y^{\mathcal{N}}-y\right\|_{L^2(\Omega)}\right), & \forall h.
\end{aligned}
\end{equation}
\end{Theorem}


\subsection{SUPG for Time-Dependent Advection-Dominated OCP($\boldsymbol{\mu}$)s} \label{time-dep-discr}
We briefly discuss the SUPG technique applied to time-dependent problems. Considering \eqref{Parametrized-Parabolic-OCP}, the main challenge comes from the fact that time derivatives should also enter into stabilization setting to ensure consistency \cite{john2011error}. However, other approaches have been proposed: in \cite{seymen2014distributed}, for instance, the time-derivative is not stabilized. Our discussion follows works inherent to Graetz-Poiseuille and Propagating Front in a Square problems without optimal control \cite{pacciarini2014stabilized,rozza2016stabilized}, where stabilization is used for time derivative, too. This adds nonsymmetric terms to the discretized state and adjoint equations. To the best of our knowledge, SUPG for Parabolic OCPs in an \textit{optimize-then-discretized} approach is still a novelty element in literature from a theoretical perspective. We refer to \cite{ganesan2012operator,giere2015supg,john2011error} for SUPG applied to generic parabolic problems.

We first discretize the equation in time, considering each discrete time as a steady-state Advection-Diffusion equation, in a \textit{space-time} approach, and then stabilized it with the SUPG.
The time interval $(0,T)$ is divided in $N_t$ sub-intervals of equal length $\Delta t := t_i - t_{i-1}$, $i \in \{1, \dots , N_t \}$. On the other hand, all terms involving time-derivative go through a time discretization equivalent to a classical implicit Euler approach \cite{akman2017streamline,hinze2008hierarchical,stoll2010all,strazzullo2020pod,strazzullo2021certified,strazzullo2022pod}. The backward Euler method is used to discretize the state equation forward in time, instead the adjoint equation is discretized backward in time using the forward Euler method, which is equivalent to the backward Euler with respect to time $T-t$, for $t \in(0, T)$ \cite{eriksson1987error, strazzullo2020pod}. The global dimension of the discrete spaces is $\mathcal{N}_{tot}=3 \cdot \mathcal{N} \cdot N_t$. We recall that $Y,U,Q$ are Hilbert Spaces and that $Y^{\mathcal{N}} \equiv (Q^{\mathcal{N}})^{*}$.

For the state equation, the stabilized term added to the form related to the time derivative of the state $\frac{\partial y}{\partial t}$ and the bilinear form $a$ is the following \cite{john2011error, pacciarini2014stabilized, rozza2016stabilized}: 
$$
s\left(y^{\mathcal{N}}(t), q^{\mathcal{N}}\right)=\sum_{K \in \mathcal{T}_{h}} \delta_{K}\left(\frac{\partial y^{\mathcal{N}}(t)}{\partial t} +\left(L_{S}+L_{SS}\right) y^{\mathcal{N}}(t), \frac{h_{K}}{|\boldsymbol{b}|} L_{S S} q^{\mathcal{N}}\right)_{K},
$$
where $y^{\mathcal{N}}(t) \in Y^{\mathcal{N}}$ for each $t \in (0,T)$ and $q^{\mathcal{N}} \in Y^{\mathcal{N}}$. Instead, the stabilized term added to the form related to the time derivative of the adjoint  $\frac{\partial p}{\partial t}$ and the bilinear form $a^{*}$ is: 
$$
s^{*}\left(z^{\mathcal{N}},p^{\mathcal{N}}(t)\right)=\sum_{K \in \mathcal{T}_{h}} \delta_{K}\left(-\frac{\partial p^{\mathcal{N}}(t)}{\partial t}+\left(L_{S}-L_{SS}\right) p^{\mathcal{N}}(t), -\frac{h_{K}}{|\boldsymbol{b}|} L_{S S} z^{\mathcal{N}}\right)_{K}.
$$

We can write the discretized state formulation using a backward Euler approach as follows:
\begin{equation}
\begin{aligned}
&\text { for each } i \in \{1, 2, \cdots, N_t\} \text {, find }  y_{i}^{\mathcal{N}} \in Y^{\mathcal{N}} \text { s.t. } \forall q^{\mathcal{N}} \in Y^{\mathcal{N}},\\
&\frac{1}{\Delta t} m_{s}\left(y_{i}^{\mathcal{N}}(\boldsymbol{\mu})-y_{i-1}^{\mathcal{N}}( \boldsymbol{\mu}), q^{\mathcal{N}} ; \boldsymbol{\mu} \right)+a_{s}\left(y_{i}^{\mathcal{N}}( \boldsymbol{\mu}), q^{\mathcal{N}} ; \boldsymbol{\mu}\right)+b_{s}\left(u^{\mathcal{N}}_i,q^{\mathcal{N}}; \boldsymbol{\mu}\right)= F_s\left(q^{\mathcal{N}}; \boldsymbol{\mu}\right),
\end{aligned}
\end{equation}
given the initial condition $y^{\mathcal{N}}_{0}$ which satisfies
\begin{equation}
\left(y^{\mathcal{N}}_{0}, q^{\mathcal{N}}\right)_{L^{2}(\Omega)}=\left(y_{0}, q^{\mathcal{N}}\right)_{L^{2}(\Omega)}, \quad \forall q^{\mathcal{N}} \in Y^{\mathcal{N}} .
\end{equation}
The stabilized term $m_s$ above is defined as:
\begin{equation}
    m_{s}\left(y^{\mathcal{N}}, q^{\mathcal{N}}; \boldsymbol{\mu}\right)=\left(y^{\mathcal{N}}, q^{\mathcal{N}}\right)_{L^{2}(\Omega)}+\sum_{K \in \mathcal{T}_{h}} \delta_{K}\left(y^{\mathcal{N}}, \frac{h_{K}}{|\boldsymbol{b}|} L_{S S} q^{\mathcal{N}}\right)_{K}
\end{equation}
and it is related to the time discretization; instead, $a_s$ and $F_s$ are defined as in the steady case. Similarly we can derive the same for the adjoint forms applying a forward Euler method:
\begin{equation}\label{eq:adjoint-discr}
\begin{aligned}
&\text { for each } i \in \{N_t-1, N_t-2,..., 1\} \text {, find } p_{i}^{\mathcal{N}} \in Y^{\mathcal{N}} \text { s.t. }\\
&\frac{1}{\Delta t} m^{*}_{s}\left(p_{i}^{\mathcal{N}}(\boldsymbol{\mu})-p_{i+1}^{\mathcal{N}}( \boldsymbol{\mu}), z^{\mathcal{N}} ; \boldsymbol{\mu} \right)+a^{*}_{s}\left(z^{\mathcal{N}}, p_{i}^{\mathcal{N}}( \boldsymbol{\mu}); \boldsymbol{\mu}\right)= - \big( y^{\mathcal{N}}_{i}-y_{d_i}, z^{\mathcal{N}}\big)_{s} \quad \forall z^{\mathcal{N}} \in Y^{\mathcal{N}}.
\end{aligned}
\end{equation}
The stabilized term $m^{*}_s$ above is defined as:
\begin{equation}
    m^{*}_{s}\left(p^{\mathcal{N}}, z^{\mathcal{N}}; \boldsymbol{\mu}\right)=\left(p^{\mathcal{N}}, z^{\mathcal{N}}\right)_{L^{2}(\Omega)}-\sum_{K \in \mathcal{T}_{h}} \delta_{K}\left(p^{\mathcal{N}}, \frac{h_{K}}{|\boldsymbol{b}|} L_{S S} z^{\mathcal{N}}\right)_{K}.
\end{equation}
Now we look at the discretization scheme. As in the steady case, $y_{i}\in Y^{\mathcal{N}}$, $u_{i}\in U^{\mathcal{N}}$ and $p_{i}\in Y^{\mathcal{N}}$, for $1 \leq i \leq N_{t}$, represent the column vectors including the coefficients of the FEM discretization for state, control and adjoint, respectively.
Therefore, we define $\boldsymbol{y}=\left[y_{1}^{T}, \ldots, y_{N_{t}}^{T}\right]^{T}$, $\boldsymbol{u}=\left[u_{1}^{T}, \ldots, u_{N_{t}}^{T}\right]^{T}$ and $\boldsymbol{p}=\left[p_{1}^{T}, \ldots, p_{N_{t}}^{T}\right]^{T}$. The vector $\boldsymbol{f}_s=\left[f_{s_{1}}^{T}, \ldots, f_{s_{N_{t}}}^{T}\right]^{T}$ indicates the components of the stabilized forcing term, $\boldsymbol{y}_{d}=$ $\left[y_{d_{1}}^{T}, \ldots, y_{d_{N_{t}}}^{T}\right]^{T}$ is the vector made of discrete time components of our desired state solution; instead, $\boldsymbol{y}_{0}=\left[y_{0}^{T}, 0^{T}, \ldots, 0^{T}\right]^{T}$ indicates the vector of initial condition for the state, where $0$ is the zero vector in $\mathbb{R}^{\mathcal{N}}$. We remark that in our test cases, we always consider $y_0 = 0$. Thus, the block matrix system is described as follows.
\begin{itemize}
    \item \textit{State equation}.
We recall that $K_s$ and $B_s$ are the matrices associated with the stabilized bilinear forms $a_s$ and $b_s$. Using the backward Euler in time, one has to solve
\begin{equation}
   M_s y_{i}+\Delta t K_s y_{i} + \Delta t B_s u_{i}=M_s y_{i-1}+f_{s_{i}} \Delta t \quad \text { for } i \in\left\{1,2, \ldots, N_{t}\right\},
\end{equation}
where $M_s$ is the stabilized mass matrix relative to the FEM discretization of $m_s$. Therefore, recalling that $y_0 = 0$, the related block matrix subsystem is
\begin{equation*}
\underbrace{\left[\begin{array}{cccc}
M_s+\Delta t K_s & 0 & & \\
-M_s & M_s+\Delta t K_s & 0 & \\
 &  \ddots & \ddots & 0 \\
 & 0 & -M_s & M_s+\Delta t K_s
\end{array}\right]}_{\mathcal{A}_s} 
\boldsymbol{y}
{+}\Delta t \underbrace{\left[\begin{array}{ccc}
B_s &0 & \\ 
0 &\ddots& 0\\ 
& 0& B_s\end{array}\right]}_{\mathcal{C}_s}
\boldsymbol{u}
= \Delta t \boldsymbol{f}_{s},
\end{equation*}
Then, everything can 
be recast in a more compact form as
\begin{equation}
\mathcal{A}_{s} \boldsymbol{y}{+}\Delta t \mathcal{C}_s \boldsymbol{u}=\Delta t \boldsymbol{f}_s.
\end{equation}

\item \textit{Gradient equation}. We recall that $B^{T}$ indicates the mass matrix related to the $b^{*}$ form and hence at every time step we have to solve the equation
\begin{equation}
    \alpha \Delta t M u_{i} {+} \Delta t B^{T} p_{i}=0, \quad \forall i \in\left\{1,2, \ldots, N_{t}\right\},
\end{equation}
which translates into the following block system:
\begin{equation*}
\Delta t \cdot \alpha  \underbrace{\left[\begin{array}{cccc}
 M & &  & \\
 & M  &  & \\
 & \ddots & \ddots & \\
 & & & M
\end{array}\right]}_{\mathcal{M}}
\left[\begin{array}{c}
u_{1} \\
u_{2} \\
\vdots \\
u_{N_{t}}
\end{array}\right] 
{+}\Delta t \underbrace{\left[\begin{array}{cccc}
B^{T} & 0 & \cdots & \\ 
& B^{T} &  & \\ 
& & \ddots & \\ 
& & & B^{T}\end{array}\right]}_{\mathcal{C^{T}}}
\left[\begin{array}{c}
p_{1} \\ 
p_{2} \\ 
\vdots \\ 
p_{N_{t}}\end{array}\right]
=\left[\begin{array}{c}
0\\ 
0\\ 
\vdots \\ 
0\end{array}\right] .
\end{equation*}
In a vector notation we have 
\begin{equation}
    \alpha \Delta t \mathcal{M} \boldsymbol{u}{+}\Delta t \mathcal{C}^{T} \boldsymbol{p}=0.
\end{equation}
\item \textit{Adjoint equation}: we have to solve {\eqref{eq:adjoint-discr}} at each time step as follows, considering $M_{s}^{T}$ the matrix formulation of $m_{s}^{*}$:
$$
M_{s}^{T} p_{i}=M_{s}^{T} p_{i+1}+\Delta t\left(-M_{s}^{T} y_{i}-K_s^{T} p_{i}+M_{s}^{T} y_{d_{i}}\right) \quad \text { for } i \in\left\{N_{t}-1, N_{t}-2, \ldots, 1\right\}.
$$
As did in previous steps, we derive the following block system:
\begin{equation*}
\underbrace{\left[\begin{array}{cccc}
M_{s}^{T}+\Delta t K_s^{T} & -M_{s}^{T} & & \\
& \ddots & \ddots & \\
& & M_{s}^{T}+\Delta t K_s^{T} & -M_{s}^{T} \\
& & & M_{s}^{T}+\Delta t K_s^{T} 
\end{array}\right]}_{\mathcal{A}^{T}_s}
\boldsymbol{p}
+\left[\begin{array}{c}
\Delta t M_{s}^{T} y_{1} \\
\vdots \\
\vdots \\
\Delta t M_{s}^{T} y_{N_{t}}
\end{array}\right]  \\
{=\left[\begin{array}{c}
\Delta t M_{s}^{T} y_{d_{1}} \\
\vdots \\
\vdots \\
\Delta t M_{s}^{T} y_{d_{N_{t}}}
\end{array}\right]}.
\end{equation*}
Then, defining $\mathcal{M}^{T}_s$ as the diagonal matrix in $\mathbb{R}^{\mathcal{N} \cdot N_{t}} \times \mathbb{R}^{\mathcal{N} \cdot N_{t}}$ which diagonal entries are $[ M_{s}^{T}, \ldots,  M_{s}^{T}]$, the adjoint system to be solved is:
$$
\Delta t \mathcal{M}^{T}_s \boldsymbol{y}+\mathcal{A}^{T}_s \boldsymbol{p}=\Delta t \mathcal{M}^{T}_s \boldsymbol{y}_{d}.
$$
\end{itemize}
 
In the end, we solve system \eqref{stab-par-block} via an one shot approach:
\begin{equation}\label{stab-par-block}
\left[\begin{array}{ccc}
\Delta t \mathcal{M}^{T}_s & 0 & \mathcal{A}^{T}_s \\
0 & \alpha \Delta t \mathcal{M} &\Delta t \mathcal{C}^{T} \\
\mathcal{A}_s &\Delta t \mathcal{C}_s & 0
\end{array}\right]\left[\begin{array}{l}
\boldsymbol{y} \\
\boldsymbol{u} \\
\boldsymbol{p}
\end{array}\right]=\left[\begin{array}{c}
\Delta t \mathcal{M}^{T}_s \boldsymbol{y}_{d} \\
0 \\
\Delta t \boldsymbol{f}_s
\end{array}\right].
\end{equation}

\section{ROMs for advection-dominated OCP($\boldsymbol \mu$)s}
\label{sec:ROMs}

{FEM simulations can be expensive in terms of computational time and memory storage}: this {issue} is obviously more evident in case of high-dimensional discrete spaces. Moreover, when we talk about parametrized PDEs, one can require to repeat simulations for several values of the parameter $\boldsymbol{\mu}$. To overcome these difficulties, we will use ROMs approach.
The basic idea of ROMs is to create a low-dimensional space, called the \textit{reduced space}, exploiting the parameter dependence of the problem at hand, such that it is a good approximation of the {initial FEM space} \cite{benner2017model,hesthaven2016certified,quarteroni2011certified,quarteroni2014reduced,quarteroni2015reduced}.
Let us consider a generic Parametrized OCPs described by the optimality conditions \eqref{Opt-lin_system}. We define the set of the parametric solutions of the optimality system with respect to the functional space $\mathbb{W}=Y \times U \times Q^{*}$ for steady $\operatorname{OCP(\boldsymbol{\mu})s}$ and  $\mathbb{W}=\mathcal{Y}_{t} \times \mathcal{U} \times \mathcal{Y}_{t}$ for the unsteady ones as
\begin{equation}
\mathbb{M}:=\{\left(y(\boldsymbol{\mu}), u(\boldsymbol{\mu}), p(\boldsymbol{\mu})\right) \ \text{solution of} \ \eqref{Opt-lin_system} \mid \boldsymbol{\mu} \in \mathcal{P}\}.
\end{equation}
The extension to space-time formulation for time-dependent problem is straightforward \cite{BALLARIN2022307, strazzullo2020pod} and requires small modifications, thus, we will exclusively refer to the steady framework.

\begin{Assumption}[Smoothness of the solution manifold]
  The continuous solution manifold $\mathbb{M}$ is smooth with respect to the parameter $\boldsymbol{\mu}  \in \mathcal{P}$.
\end{Assumption} 

{In order to execute computational simulations, we want to find a suitable discretization $\mathbb{M}^{\mathcal{N}}$ of $\mathbb{M}$: as the numerical solution of \eqref{Opt-lin_system} can present oscillations that do not approximate correctly the true solution, the manifold $\mathbb{M}^{\mathcal{N}}$ will be described by solutions of the stabilized system \eqref{stab-block-system}.} 
Let $\mathbb{W}^{\mathcal{N}} \subset \mathbb{W}$ be our FEM approximation of the continuous space $\mathbb{W}$, we call $\mathbb{W}^{\mathcal{N}}:= Y^{\mathcal{N}} \times U^{\mathcal{N}} \times \big(Q^{\mathcal{N}}\big)^{*}$ the high-fidelity space. Then, for stabilized problems we define the discrete parametric solution manifold as
\begin{equation}\label{hf-sampled}
\mathbb{M}^{\mathcal{N}}:=\left\{\left(y^{\mathcal{N}}(\boldsymbol{\mu}), u^{\mathcal{N}}(\boldsymbol{\mu}), p^{\mathcal{N}}(\boldsymbol{\mu})\right) \text{ FEM solution of the } \eqref{stab-block-system} \mid \boldsymbol{\mu} \in \mathcal{P}\right\}.
\end{equation}
Starting from $\mathbb{M}^{\mathcal{N}}$, ROM techniques create a reduced space of low dimension $N$ denoted with $\mathbb{W}^{N}$, via a linear combination of \textit{snapshots}, i.e.\ high-fidelity evaluations of the optimal solution $\left(y^{\mathcal{N}}(\boldsymbol{\mu}), u^{\mathcal{N}}(\boldsymbol{\mu}), p^{\mathcal{N}}(\boldsymbol{\mu})\right)$ computed in properly chosen parameters values $\boldsymbol{\mu}$. Obviously we have that $\mathbb{W}^{N} \subset \mathbb{W}^{\mathcal{N}}$ and we denote $\mathbb{W}^{N} =  Y^{N} \times U^{N} \times \big(Q^{N}\big)^{*}$. Here, $ Y^{N}$, $U^{N}$ and $(Q^{N})^{*}$ are the \textit{reduced spaces} for the state, the control and the adjoint variables, respectively. The snapshots are collected by a POD algorithm using a \textit{partitioned approach}. This strategy is followed due to good results shown in literature \cite{karcher2018certified,negri2015reduced,strazzullo2018model}. After having built these reduced function spaces, a standard Galerkin projection is performed onto these ones \cite{atkinson2005theoretical, quarteroni2009numerical, quarteroni2008numerical}.

\subsection{Offline-Online Procedure for ROMs}
\label{sec: rom procedures}
ROM procedure is divided in two stages:
\begin{itemize}
    \item \textit{offline} phase: here the \textit{snapshots} {are collected} by solving the{stabilized} high-fidelity system
\eqref{stab-block-system}. Secondly, the low-dimensional bases are created and hence all reduced spaces $Y^{N}$, $U^{N}$ and $(Q^{N})^{*}$ are built and stored, too. Moreover, all the $\boldsymbol{\mu}$-independent quantities are assembled and stored. It is potentially an expensive phase, which depends on $\mathcal{N}$.
    \item \textit{online} phase: here a parameter $\boldsymbol{\mu}$ is chosen and all the previous stored $\boldsymbol{\mu}$-independent quantities are {combined} with the just-computed $\boldsymbol{\mu}$-dependent ones to build the reduced block matrix system based on a Galerkin projection. To be convenient, this phase should be $\mathcal{N}$-independent. Whereas in the offline phase stabilization is present due to stabilized snapshots, for the online phase this might not be necessary. Therefore, we have two possibilities: if stabilization is performed also here, we talk about \textit{Online-Offline stabilization}, otherwise we denote the setting as \textit{Only-Offline stabilization}.
\end{itemize}
A sufficient condition to perform the online phase independently of $\mathcal{N}$ is to admit the separation of the variables depending on $\boldsymbol{\mu}$ and the solution $(y,u,p)$ in the \textit{affine decomposition} \cite{hesthaven2016certified}. 
\begin{Assumption}\label{affine-ass}
We require that all the forms in{the stabilized system} \eqref{stab-block-system} are affine in $\boldsymbol{\mu} \in \mathcal{P}$.
\end{Assumption}
In Section \ref{sec:pod} we describe the POD algorithm used in the offline phase. Now, we illustrate the explicit expression of the reduced solutions. In terms of their basis functions, we define
\begin{equation}\label{basis-red}
\hspace{-0.1cm}Y^{N} =\operatorname{span}\left\{\eta_{n}^{y}, \ n=1,  \ldots, N\right\},  \ 
U^{N} =\operatorname{span}\left\{\eta_{n}^{u}, \ n=1, \ldots, N\right\}, \
(Q^{N})^{*} =\operatorname{span}\left\{\eta_{n}^{p},\  n=1, \ldots, N\right\},
\end{equation}
the reduced state, the reduced control and the reduced adjoint space, respectively. After having built them, we consider an enriched space for state and adjoint variables. 
Therefore, let us denote with $\left\{\tau_{n}\right\}_{n=1}^{2 N}=\left\{\eta_{n}^{y}\right\}_{n=1}^{N} \cup\left\{\eta_{n}^{p}\right\}_{n=1}^{N}$ the basis functions for the space $Z^{N}$, with $Z^N \equiv Y^N\equiv (Q^{N})^{*}$, then we have
$
Z^{N}=\operatorname{span}\left\{\tau_{n}, n=1, \ldots, 2 N\right\}
$ \cite{doi:10.1137/090760453,doi:10.1137/110854084,karcher2018certified,kunisch2008proper, negri2013reduced,negri2015reduced}.

\subsection{Proper Orthogonal Decomposition}
\label{sec:pod}
In this Section we briefly describe the Proper Orthogonal Decomposition (POD) Galerkin algorithm \cite{BALLARIN2022307,hesthaven2016certified,strazzullo2018model,strazzullo2020pod} for the construction of a discrete solution manifold and the relative reduced spaces. Since in the unsteady case we use a space-time structure, this procedure can be described making no distinction between time-dependency and steadiness. 
Firstly, we make a sampling of $\mathcal{P}$ by choosing $N_{\text{train}}$ of its elements. Therefore, let us define the set of the train samples as $\mathcal{P}_{N_{\text{train}}}$: we have that obviously $\mathcal{P}_{N_{\text{train}}} \subset \mathcal{P}$ and the cardinality is $|\mathcal{P}_{N_{\text{train}}} | = N_{\text{train}}$. The set $\mathcal{P}_{N_{\text{train}}}$ is denoted as the \textit{training set}. We should pursue that $N_{\text{train}}$ is
large enough so as to ensure that $\mathcal{P}_{N_{\text{train}}}$ is a good “approximation” of the parameter space $\mathcal{P}$. $\mathcal{P}_{N_{\text{train}}}$ is built through a Monte-Carlo sampling method with respect to a uniform
density with support equal to $\mathcal{P}$.
  
Starting from the sampling, the POD algorithm manipulates $N_{\text{train}}$ snapshots for the state, the adjoint and the control variables:
\begin{equation}\label{snapshots}
\left\{\left(y^{\mathcal{N}}(\boldsymbol{\mu}_{j}), u^{\mathcal{N}}(\boldsymbol{\mu}_{j}), p^{\mathcal{N}}(\boldsymbol{\mu}_{j})\right)\right\}_{j=1}^{N_{\text{train}}} \text{ with } \boldsymbol{\mu}_{j} \in \mathcal{P}_{N_{\text{train}}}.
\end{equation}
After this step, a compressing stage is performed: from \eqref{snapshots} we build $N$ basis functions by only considering the most important parametric information and throwing away the redundant ones, with $N \leq N_{\text{train}}$. 

Let us consider $W^{\mathcal N}$ a generic high-fidelity space.
Given parametric solutions $w^{\mathcal N}(\boldsymbol \mu_j) \in W^{\mathcal N}$, for $1 \leq j \leq N_{train}$, the POD provide a low-dimensional space $W_N \subset W^{\mathcal N}$ minimizing the following quantity:
\begin{equation*}
 \sqrt{\frac{1}{N_{\text{train}}} \sum_{\boldsymbol{\mu}_{j} \in \mathcal{P}_{N_{\text{train}}}} \min _{\bar{w} \in {W}^{N}}\left\| w^{\mathcal{N}}\left(\boldsymbol{\mu}_{j}\right)- \bar{w}\right\|_{W}^{2}}.
\end{equation*}
Firstly, we collect a set of ordered parameters $\boldsymbol{\mu}_{1}, \ldots, \boldsymbol{\mu}_{N_{\text{train}}} \in \mathcal{P}_{N_{\text{train}}}$ and compute the ordered snapshots $w^{\mathcal{N}}\left(\boldsymbol{\mu}_{1}\right), \ldots, w^{\mathcal{N}}\left(\boldsymbol{\mu}_{N_{\text{train}}}\right)$. Let us define $\mathbf{C} \in \mathbb{R}^{N_{\text{train}} \times N_{\text{train}}}$ as the correlation matrix of the snapshots as follows:
\begin{equation}\label{correlation-matrix}
\mathbf{C}_{ij} :=\frac{1}{N_{\text{train}}}\left(w^{\mathcal{N}}\left(\boldsymbol{\mu}_{i}\right), w^{\mathcal{N}}\left(\boldsymbol{\mu}_{j}\right)\right)_{W}, \quad 1 \leq i,j \leq N_{\text{train}}.
\end{equation}
The next step is to find the pair eigenvalue-eigenvector $\left(\lambda_{n}, e_{n}\right)$ of the following problem:
$$
\mathbf{C} e_{n} =\lambda_{n} e_{n}, \mbox{ with } \| e_{n}\| = 1, \mbox{ for } 1 \leq n \leq N_{\text{train}}.
$$
For the sake of simplicity, we organise the eigenvalues $\lambda_{1}, \ldots, \lambda_{N_{\text{train}}}$ in decreasing order. Consider the first $N$ ones, specifically $\lambda_{1}, \ldots, \lambda_{N}$ together with the related eigenvectors $e_{1}, \ldots, e_{N}$. We refer to the  $k$-th component of the state eigenvector $e_{n} \in \mathbb{R}^{N_{\text{train}}}$ with the notation $\left(e_{n}\right)_{k}$. After having finished this step, the basis functions $\eta_{n}$ for reduced space are built through the following formula:
\begin{equation}\label{eigen-const}
\eta_{n} =\frac{1}{\sqrt{\lambda_{n}}} \sum_{k=1}^{N_{\text{train}}}\left(e_{n}\right)_{k} w^{\mathcal{N}}\left(\boldsymbol{\mu}_{k}\right), \quad 1 \leq n \leq N.
\end{equation}
The POD algorithm is summarized in Algorithm \ref{alg:POD}. For the OCP problem, we employed a \textit{partitioned approach}: we perform the POD algorithm separately for all three variables, compressing $Y^{\mathcal N}$ in $Y^N$, $U^{\mathcal N}$ in $U^N$ and $(Q^{\mathcal N})^{*}$ in $(Q^{N})^{*}$. Then, the \emph{aggregated space} technique is applied and we combine the reduced state space $Y^N$ and the adjoint space $(Q^{N})^{*}$ in the space $Z^N$. The latter space will be used to describe both state and adjoint variables.

As both $N$ and $N_{\text{train}}$ can be chosen by us, we should find sharp criteria in order to decide them. A possibility can be to set them in based on a study of the eigenvalues, using the estimate \cite{hesthaven2016certified, quarteroni2015reduced, rozzaweighted}:
\begin{equation}\label{eigenvalues}
    \sqrt{\frac{1}{N_{\text{train}}} \sum_{k=1}^{N_{\text{train}}}\left\|w^{\mathcal{N}}\left(\boldsymbol{\mu}_{k}\right)-\Pi^{N}\left(w^{\mathcal{N}}\left(\boldsymbol{\mu}_{k}\right)\right)\right\|_{{W}}^{2}}=\sqrt{\sum_{k=N+1}^{N_{\text{train}}} \lambda_{k}},
\end{equation}
where $\Pi^{N}: W \rightarrow W^{N}$ is a Galerkin projector of functions from $W$ onto $W^{N}$. 
The second member of equation \eqref{eigenvalues} can be a measure of how well the FEM space is approximated by $N$ reduced basis over the chosen training set of cardinality $N_{\text{train}}$.

\begin{algorithm}
\caption{POD algorithm}\label{alg:POD} 
\textbf{Input}: parameter domain $\mathcal{P}$, FEM spaces $W^{\mathcal{N}}$ and $N_{\text{train}}$. 

\textbf{Output:} reduced spaces $W^{N}$. 

Starting from the high-fidelity spaces $W^{\mathcal{N}}$

\begin{algorithmic}[1]
\State Sample $\mathcal{P}_{\text{train}} \subset \mathcal{P}$;
\ForAll {$\boldsymbol{\mu} \in \mathcal{P}_{\text{train}}$} 
\State Collect high-fidelity snapshots $w^{\mathcal{N}}\left(\boldsymbol{\mu}_{i}\right) \in W^{\mathcal N}$;
\EndFor
\State Assemble the matrix $\mathbf{C}_{ij} :=\frac{1}{N_{\text{train}}}\left(w^{\mathcal{N}}\left(\boldsymbol{\mu}_{i}\right), w^{\mathcal{N}}\left(\boldsymbol{\mu}_{j}\right)\right)_{W}, \ 1 \leq i,j \leq N_{\text{train}}$.
\State Compute the eigenvalues $\lambda_{1}, \ldots, \lambda_{N_{\text{train}}}$ and the corresponding orthonormalized eigenvectors $e_{1}, \ldots, e_{N_{\text{train}}}$ of $\mathbf{C}$;
\State {After having chosen $N$ according to a certain criterion}, define $W^{N}= \operatorname{span}\left\{\eta_{n}, n=1, \ldots, N\right\}$, where $\eta_{n} =\frac{1}{\sqrt{\lambda_{n}}} \sum_{k=1}^{N_{\text{train}}}\left(e_{n}\right)_{k} w^{\mathcal{N}}\left(\boldsymbol{\mu}_{k}\right)$.
\end{algorithmic}
\end{algorithm}
}

\begin{Remark}[Time-dependent problems]
  When we are dealing with time-dependent OCPs, the time instances are not separated in the POD procedure. Therefore, the space-time problem is studied as a steady one and each snapshot carries {all the time instances}.
\end{Remark}

\subsection{Only-Offline vs Online-Offline stabilization}
We briefly illustrate what practically means to perform the two different ROM procedures exposed in Section \ref{sec: rom procedures}. The two techniques are shown only for steady problems, but the extension to parabolic ones is straightforward.

The offline phase is common for both procedures. The snapshots are collected by solving the stabilized system \eqref{supg-system}. Then, the FEM spaces $Y^N$, $U^N$, and $(Q^N)^{*}$ are computed according to the POD procedure in Algorithm \ref{alg:POD} and, straight after, the enriched space $Z^N$ is also built. What makes different the two strategies is the online phase.
The Only-Offline method aims to find a solution of the following reduced order model system:
\begin{equation}\label{rom onlyoff-system}
    \begin{aligned}
    a^{*}\left(z^N, p^N\right)+\big( y^N-y_d, z^N\big)=0 , \qquad &\forall z^N \in Z^N,\\
   b^{*}\big(v^N, p^N\big)+\alpha n\big(u^N, v^N\big)=0 , \qquad & \forall v^N \in U^N, \\
         a\left(y^N, q^N\right)+b\left(u^N,q^N\right)=F(q^N) , \qquad & \forall q^N \in Z^N,
    \end{aligned}
\end{equation}
which means considering system \eqref{Opt-lin_system} without any stabilization in the online phase, in contrast to the offline one. Instead, in the Offline-Online procedure we aim to solve the stabilized system in the online phase, too:
\begin{equation}\label{rom onoff-system}
    \begin{aligned}
   a_{s}^{*}\left(z^N, p^N\right)+\big( y^N-y_d, z^N\big)_{s}=0 , \qquad &\forall z^N \in Z^N,\\
   b^{*}\big(v^N, p^N\big)+\alpha n\big(u^N, v^N\big)=0 , \qquad & \forall v^N \in U^N, \\
  a_{s}\left(y^N, q^N\right)+b_{s}\left(u^N,q^N\right)=F_{s}(q^N) , \qquad & \forall q^N \in Z^N.
    \end{aligned}
\end{equation}
It is worth mentioning that both state and adjoint reduced solutions are sought in the aggregated space $Z^N$.

\section{Numerical Results}
\label{sec:results}

In this Section we propose simulations regarding the Graetz-Poiseuille and the Propagating Front in a Square problems.
Regarding the steady case, the numerical experiments are coded through the RBniCS library \cite{RBniCS}; instead, the unsteady ones are implemented employing both RBniCS and multiphenics \cite{multiphenics} libraries. They are python-based libraries, built on FEniCS \cite{logg2012automated}.
When we perform the Online-Offline stabilization procedure, we will always use the same stabilization parameter $\delta_K$ of the high-fidelity approximation also at the reduced level, both in steady and unsteady cases.

We will illustrate an analysis over \emph{relative errors} between the FEM and the reduced solutions for all three variables, defined as
\begin{equation} \label{rel-error}
e_{y, N}(\boldsymbol{\mu}):=\frac{\left\|y^{\mathcal{N}}(\boldsymbol{\mu})-y^{N}(\boldsymbol{\mu})\right\|_{Y}}{\left\|y^{\mathcal{N}}(\boldsymbol{\mu})\right\|_{Y}}, \
e_{u, N}(\boldsymbol{\mu}):=\frac{\left\|u^{\mathcal{N}}(\boldsymbol{\mu})-u^{N}(\boldsymbol{\mu})\right\|_{U}}{\left\|u^{\mathcal{\mathcal{N}}}(\boldsymbol{\mu})\right\|_{U}}, \
e_{p, N}(\boldsymbol{\mu}):=\frac{\left\|p^{\mathcal{N}}(\boldsymbol{\mu})-p^{N}(\boldsymbol{\mu})\right\|_{Q^{*}}}{\left\|p^{\mathcal{N}}(\boldsymbol{\mu})\right\|_{Q^{*}}},
\end{equation}
for the state, the control and the adjoint, respectively.
As we are dealing with parametrized OCPs, we will evaluate a simple average of \eqref{rel-error} for $\boldsymbol{\mu}$ {uniformly distributed} in a testing set $\mathcal{P}_{\text{test}} \subseteq \mathcal{P}$  of size $N_{\text{test}}$ for every dimension $N=1, \ldots, N_{\max}$ of the reduced space obtained by our POD procedure.
More precisely, we will plot the base-$10$ logarithm of the average of \eqref{rel-error}. For parabolic problems, we will consider the sum of the errors with respect to each discretized instant of time $t$.

Regarding the efficiency of ROMs, we use the \textit{speedup-index} to compare the computational cost of the FEM solution with that of the reduced one. This quantity is defined as:
\begin{equation}\label{speed-up}
\text{speedup-index}=\frac{\text{ computational time of the high-fidelity solution }}{\text{ computational time of the reduced solution }},
\end{equation}
which will be computed for each $\boldsymbol{\mu}$ in the testing set with respect to the dimension $N$ of the reduced spaces. As made with the relative error, we will consider the sample average of this quantity with respect to $N$; however, for the sake of completeness, we will add its minimum and maximum value computed through the testing set. For each test case, we will use the same $\mathcal{P}_{\text{test}}$ to compute relative errors and the speedup-index. The steady results are obtained with $16$GB of RAM and Intel Core i7-7500U Dual Core, $2.7$GHz for the CPU; instead, the FEM and ROM parabolic simulations are run with $16$GB of RAM and Intel Core i$7-7700$ Quad Core, $3.60$GHz for the CPU.{For all simulations, we solve our algebraic system derived from the discretization procedure relying on the direct \textit{MUMPS} solver present in FEniCS \cite{amestoy2001fully,amestoy2019performance}.}
{\begin{Remark}
We stress that the main goal of this numerical investigation was to analyze the advantages of Online-Offline stabilization with respect to the Only-Offline ones. In the simulations, the variables seem not to completely reach the steady state behaviour. We believe that this is due, on one side, to the time dependency of the control variable, and, on the other, to the small final time $T$ values we chose for the simulations. The choice of $T$ was influenced by the computational effort needed for the simulation.
\end{Remark}}
\subsection{Numerical Experiments for the Graetz-Poiseuille Problem}
\label{sec:c graetz}

\begin{figure}[!h]
\centering
     \includegraphics[scale=0.23]{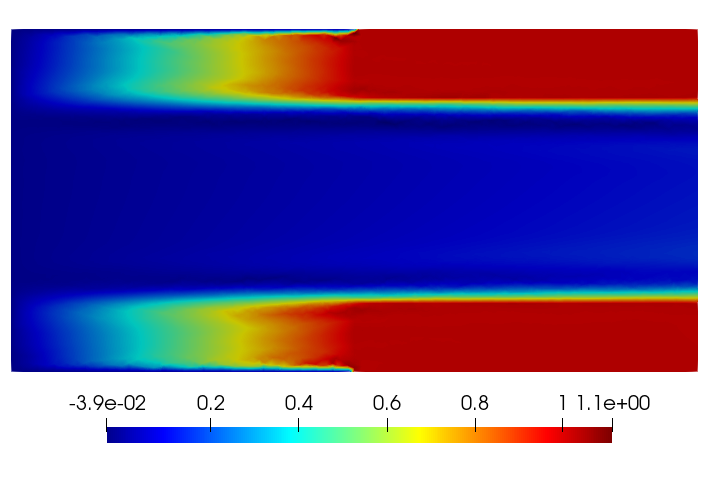} 
      \includegraphics[scale=0.225]{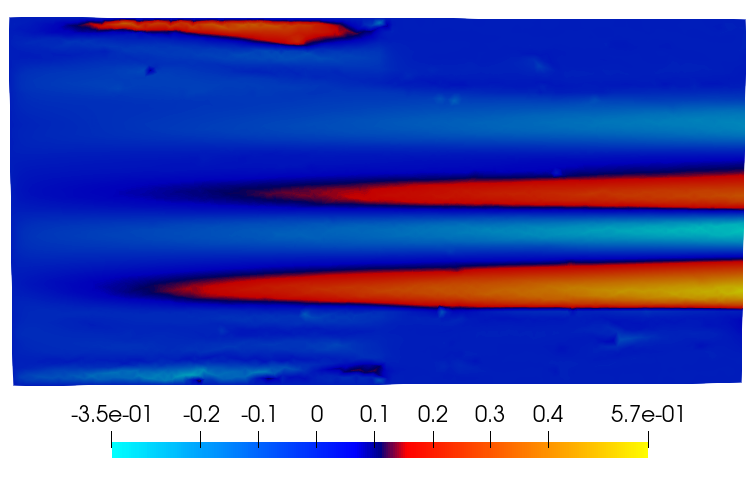} 
     \includegraphics[scale=0.23]{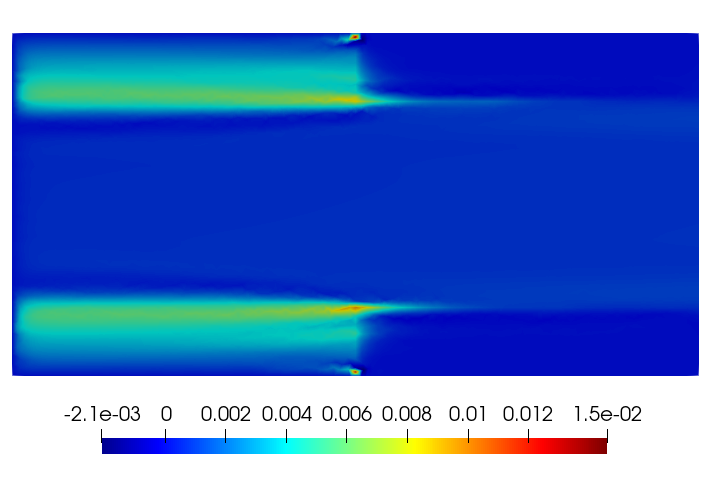} 
     \includegraphics[scale=0.23]{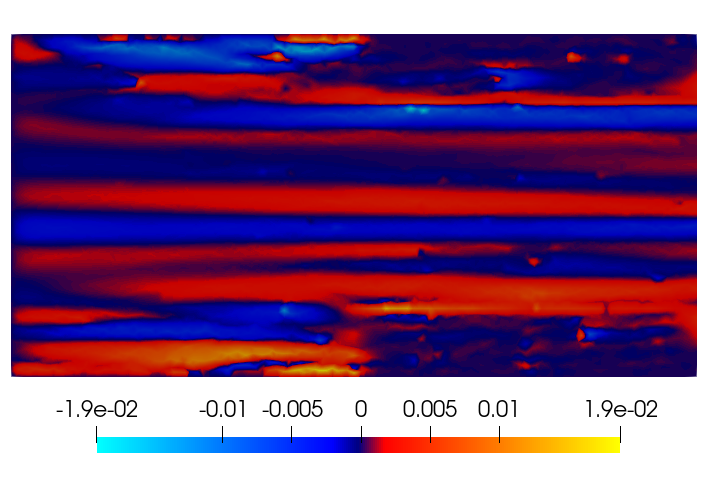} 
     \caption{{(\underline{Top}) Only-Offline stabilized state (left) and its error with respect to the FEM solution (right); (\underline{Bottom}) Only-Offline stabilized adjoint (left) and its error with respect to the FEM solution (right); for the Graetz-Poiseuille Problem; $N=1$, $\mathcal{P} = \big[10^4,10^6\big]$, $\mu_1=10^5$, $h=0.029$, $\delta_K =1.0$, $\alpha=0.01$.}}
     \label{fig:onoff-graetz-yp1}
\end{figure}

\begin{figure}[!h]
\centering
     \includegraphics[scale=0.23]{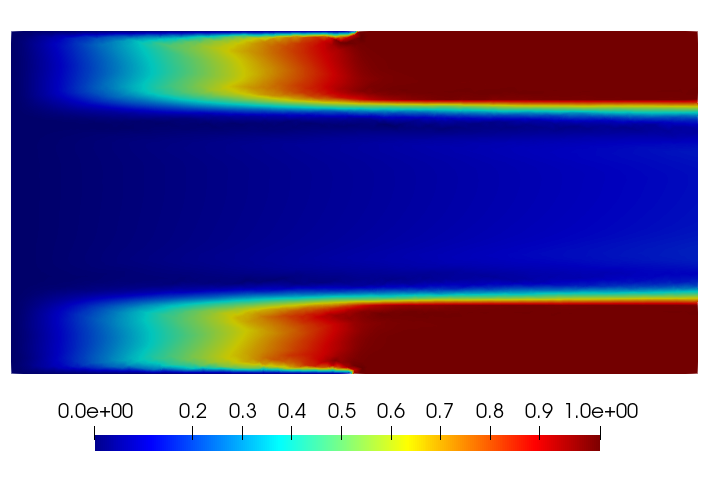}
     \includegraphics[scale=0.22]{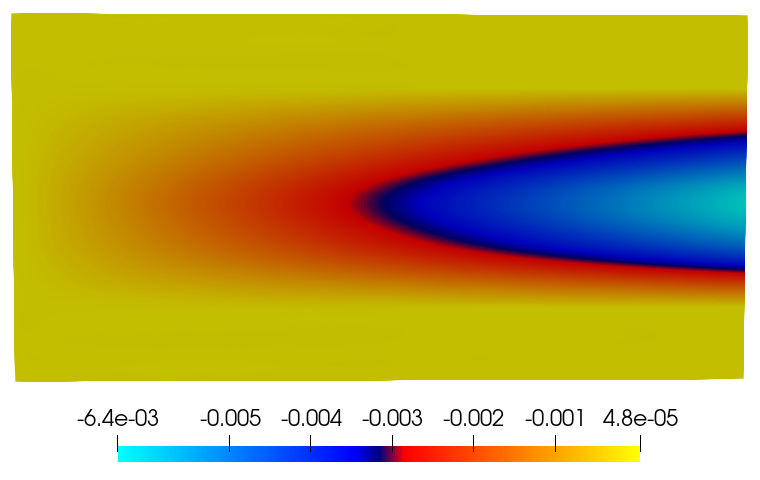}
     \includegraphics[scale=0.23]{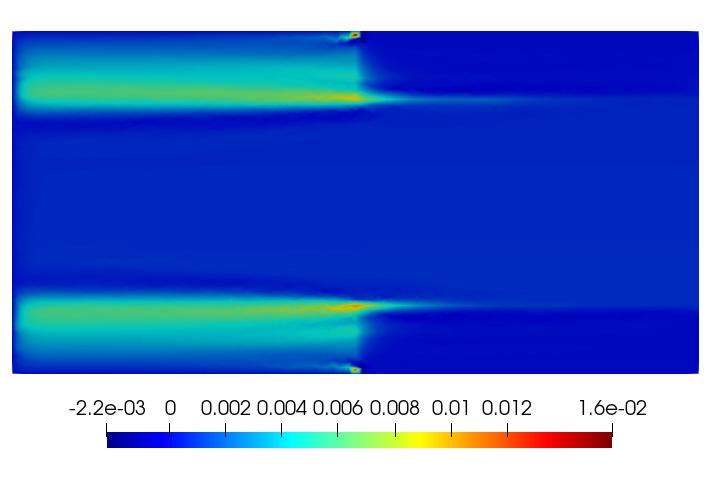}
     \includegraphics[scale=0.224]{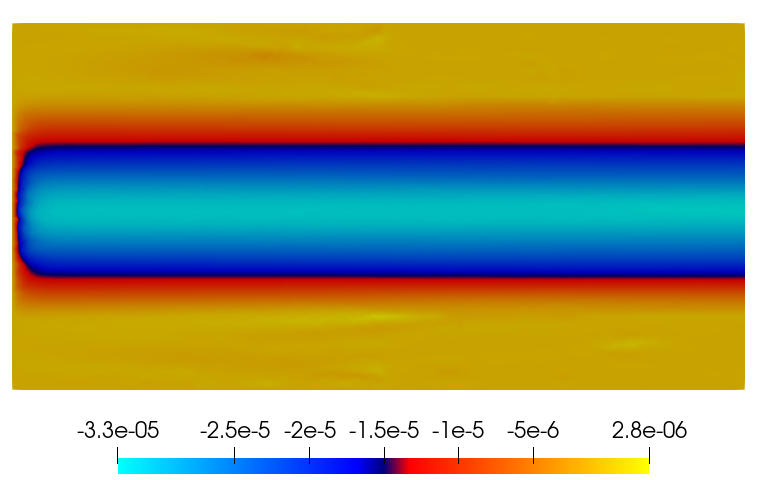}
     \caption{{(\underline{Top}) Online-Offline stabilized state (left) and its error with respect to the FEM solution (right); (\underline{Bottom}) Online-Offline stabilized adjoint (left) and its error with respect to the FEM solution (right); for the Graetz-Poiseuille Problem; $N=6$, $\mathcal{P} = \big[10^4,10^6\big]$, $\mu_1=10^5$, $h=0.029$, $\delta_K =1.0$, $\alpha=0.01$.}}
     \label{fig:onoff-graetz-yp2}
\end{figure}

\begin{figure}[!h]
\centering
     \includegraphics[scale=0.154]{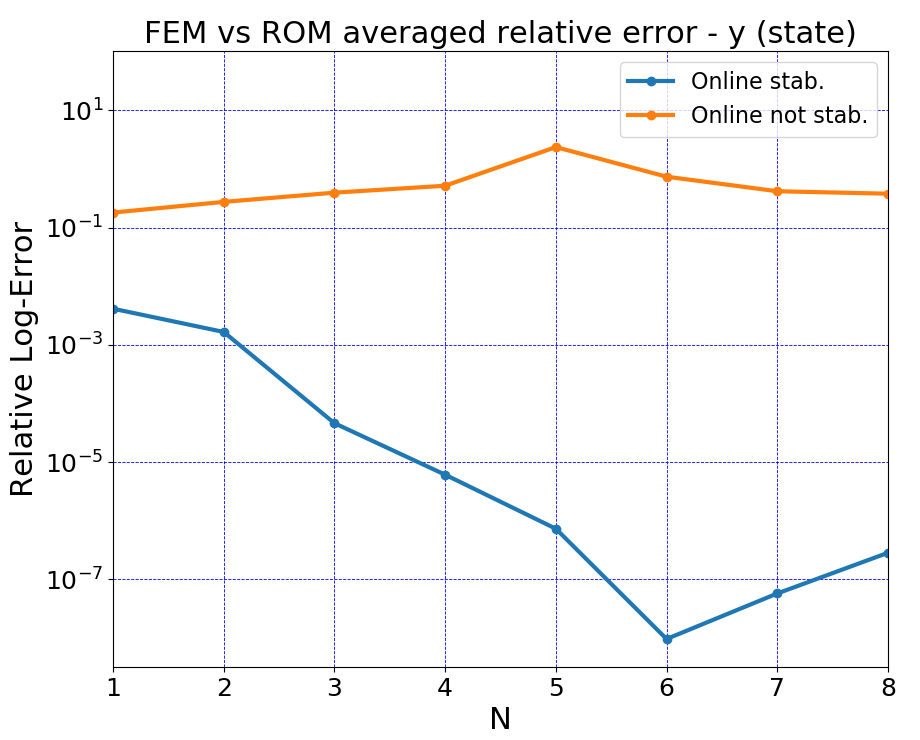} 
     \includegraphics[scale=0.154]{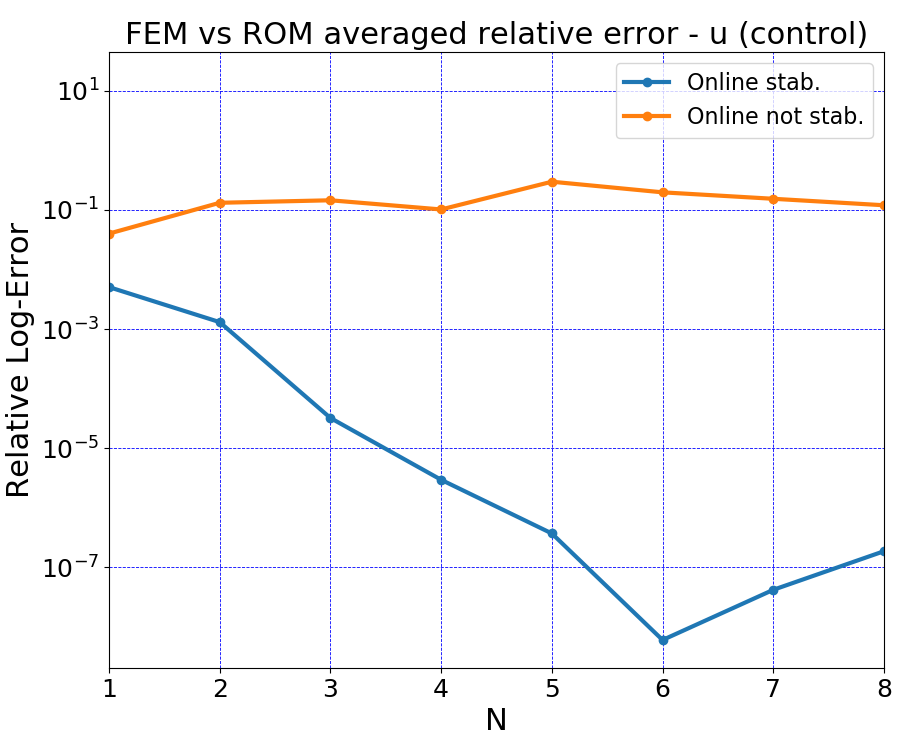}
     \includegraphics[scale=0.154]{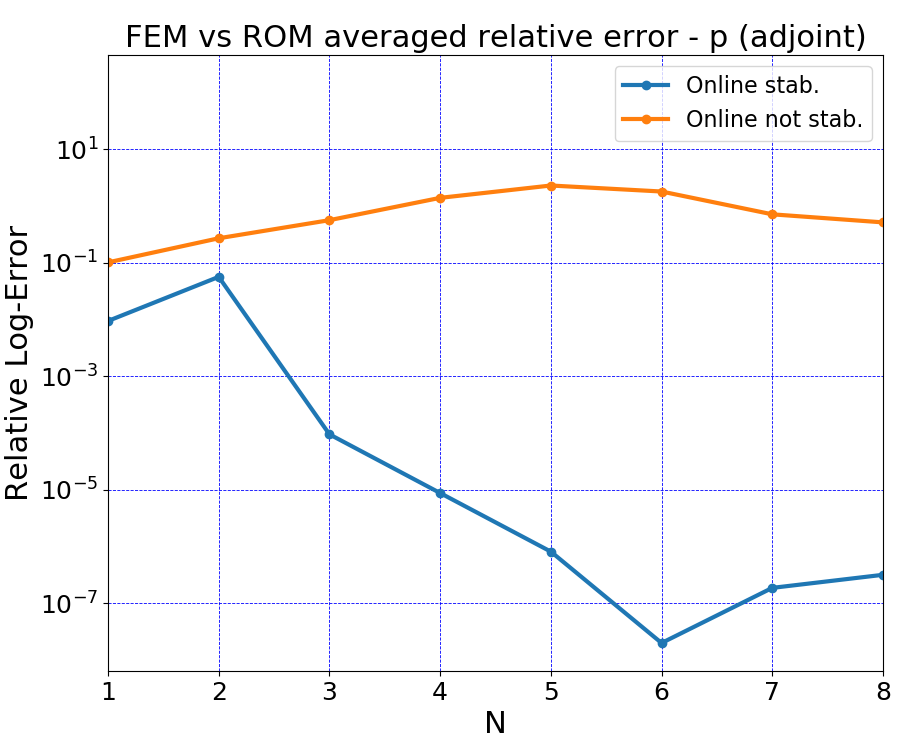}
     \caption{Relative errors between FEM and reduced solution for state (left), control (center) and adjoint (right), for Online-Offline and Only-Offline stabilization, $\alpha=0.01$, $N_{\text{test}}=100$, $h=0.029$, $\mathcal{P}=\big[10^4,10^6\big]$. Graetz-Poiseuille Problem.}
     \label{fig:error-graetz-yup}
\end{figure}

The Graetz-Poiseuille problem concerns the heat conduction in a straight duct, whose walls can be characterized by heat exchange or maintained at a certain fixed temperature. This example is very well-known in the numerical Advection-Dominated literature \cite{gelsomino2011comparison,pacciarini2014stabilized, rozza2009reduced, torlo2018stabilized}. 
We start by presenting the stationary case. We apply a distributed control in the whole domain and the parameter $\boldsymbol{\mu} = \mu_1>0$ is a physical component and characterizes the diffusion term.
\begin{figure}[h!]
   \centering
        \begin{tikzpicture}[scale=3.0]
                 \draw[color=DEblue!100, fill=DEblue!10] (1,1) -- (2, 1) -- (2,0.8) -- (1,0.8) --(1,1)node[midway, left, scale=1.2]{};	
                 \draw[color=DEblue!100, fill=DEblue!10] (1,0) -- (2, 0) -- (2,0.2) -- (1,0.2) --(1,0)node[midway, left, scale=1.2]{};
                 \draw[black] (1.5,0.9) node[scale=1.]{{$\Omega_{obs}$}};
                 \draw[black] (1.5,0.1) node[scale=1.]{{$\Omega_{obs}$}};
                 \draw[black] (1,0.5) node[scale=1.5]{{$\Omega$}};
                  \draw[blue] (0,0) -- (1,0) node[midway, below, scale=1.2]{$\Gamma_{1}$};
                 \draw[red] (1,0) -- (2,0) node[midway, below, scale=1.2]{$\Gamma_{2}$};
                 \draw[dashed,red] (2,0) -- (2,1) node[midway, right, scale=1.2]{$\Gamma_{3}$};
                 \draw[red] (2,1) -- (1,1) node[midway, above, scale=1.2]{$\Gamma_{4}$};
                 \draw[blue] (1,1) -- (0,1) node[midway, above, scale=1.2]{$\Gamma_{5}$};
                 \draw[blue] (0,1) -- (0,0) node[midway, left, scale=1.2]{$\Gamma_{6}$};
                 \filldraw[black] (0,0) circle (0.3pt) node[left]{(0,0)};
                 \filldraw[black] (1,0) circle (0.3pt) node[below]{(1,0)};
                 \filldraw[black] (2,0) circle (0.3pt) node[right]{(2,0)};
                 \filldraw[black] (2,0.2) circle (0.3pt) node[right]{(2,0.2)};
                 \filldraw[black] (2,0.8) circle (0.3pt) node[right]{(2,0.8)};
                 \filldraw[black] (2,1) circle (0.3pt) node[right]{(2,1)};
                 \filldraw[black] (1,1) circle (0.3pt) node[above]{(1,1)};
                 \filldraw[black] (0,1) circle (0.3pt) node[left]{(0,1)};
        \end{tikzpicture}
        \caption{Geometry of the Graetz-Poiseuille Problem.}
       \label{fig:Geometry-Graetz}
\end{figure}
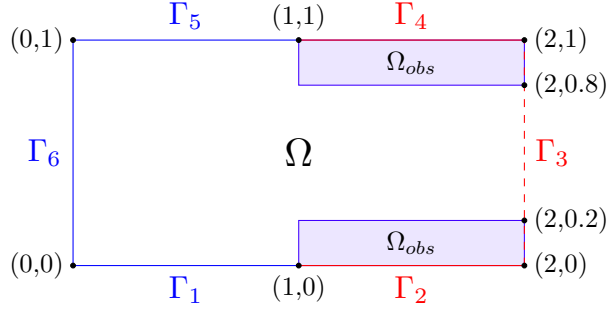
{The spatial coordinates of the system are denoted with $(x_0, x_1)$.}
{The boundary of $\Omega$ is $\Gamma$.}
We consider Dirichlet boundary conditions (BC) on sides $\Gamma_{1}:=[0,1]\times \{0\}$, $\Gamma_{5}:= [0,1]\times \{1\}$, $\Gamma_{6}:=\{0\}\times [0,1]$ {by imposing} $y=0$ and $\Gamma_{2}:=[1,2]\times \{0\}$ and $\Gamma_{4}:=[1,2]\times \{1\}$ {by imposing} $y=1$, referring to Figure \ref{fig:Geometry-Graetz}. We deal with homogeneous Neumann conditions on $\Gamma_{3}:=\{2\} \times [0,1]$. The classic formulation of the problem is:
\begin{equation} \label{graetz-system}
    \begin{cases}
      \displaystyle -\frac{1}{{\mu_{1}}} \Delta y(\boldsymbol{\mu})+4 x_1(1-x_1) \partial_{x_0} y(\boldsymbol{\mu})=u, & \text { in } \Omega, \\
y(\boldsymbol{\mu})=0, & \text { on } \Gamma_{1} \cup \Gamma_{5} \cup \Gamma_{6}, \\
y(\boldsymbol{\mu})=1, & \text { on } \Gamma_{2} \cup \Gamma_{ 4}, \\
\displaystyle \frac{\partial y(\boldsymbol{\mu})}{\partial \nu}=0, & \text { on } \Gamma_{3}.
    \end{cases}
\end{equation}
Now we want to derive the optimality system. $\Omega_{obs}:= [1,2] \times [0.8,1] \cup [1,2] \times [0,0.2]$ as illustrated in Figure \ref{fig:Geometry-Graetz}. {In this case, the state belongs to the space:
$$\tilde{Y}:= \big\{v \in H^{1} \big(\Omega\big) \text{ s.t. it satisfies the Dirichlet } \mathrm{BCs} \text{ in }  (\ref{graetz-system}) \big\} .$$
}
{
For the sake of practice, it is better to introduce a lifting function $R_{y} \in H^{1}(\Omega)$, such that it fulfills the $\mathrm{BC} \text{ in }  (\ref{graetz-system})$. Therefore we define the variable $\bar{y} := y - R_y$, with $\bar{y} \in Y$, where
$$Y:= \big\{ v \in H^{1}_{0}\big(\Omega\big) 
\text{ s.t. } 
v=0 \text{ on } \Gamma \setminus \Gamma_{3} \big\}.$$ Nevertheless, without loss of generality, we will denote the new variable $\bar{y}$ with $y$ and we settle $U := \mathrm{L}^{2}(\Omega)$ and $Q := Y^{*}$. } Therefore, the adjoint variable $p$ is null on $ \Gamma \setminus \Gamma_{3}$.  The mathematical formulation is described as follows (we omitted the dependence from $\boldsymbol{\mu}$). Fixed $\alpha>0$, find the pair $(y,u) \in Y \times U$ that realizes
\begin{equation}
\label{Jgraetz}
\displaystyle
    \min\limits_{(y,u) \in Y \times U} J(y,u) = \frac{1}{2} \int\limits_{\Omega_{obs}} 
    \big(y -y_d\big)^2 \ \mathrm{dx} +\frac{\alpha}{2} \int\limits_{\Omega}  u^2 \ \mathrm{dx} \quad \text{such that} \ e\left(y, u, p ; \boldsymbol{\mu}\right)=0,
\end{equation}
where $
    e\left(y, u, p ; \boldsymbol{\mu}\right) := a\left(y, p ; \boldsymbol{\mu}\right) + b\left(u, p ; \boldsymbol{\mu}\right)-\langle p, f(\boldsymbol{\mu}) \rangle_{Y^{*}Y}.$
As explained in Sections \ref{sec:problem} and \ref{sec:truth_discretization}, we follow a Lagrangian approach and we use SUPG stabilization in a \emph{optimize-then-discretize} framework.
We exploit $\mathbb{P}^{1}$-FEM approximation for the state, control and adjoint spaces. Here the stabilized forms $a_s$ and $a^{*}_s$ are, respectively:
\begin{equation*}
  \begin{aligned}
a_{s}\left(y^{\mathcal{N}}, q^{\mathcal{N}} ; \boldsymbol{\mu}\right)&:=a\left(y^{\mathcal{N}}, q^{\mathcal{N}} ; \boldsymbol{\mu}\right)+\sum_{K \in \mathcal{T}_{h}} \delta_{K} \int_{K}\left(4 x_1(1-x_1) \partial_{x_0} y^{\mathcal{N}}\right)\left(h_{K} \partial_{x_0} q^{\mathcal{N}}\right), \quad y^{\mathcal{N}}, q^{\mathcal{N}} \in Y^{\mathcal{N}},\\
   a^{*}_{s}\left(z^{\mathcal{N}}, p^{\mathcal{N}} ; \boldsymbol{\mu}\right)&:=a^{*}\left(z^{\mathcal{N}}, p^{\mathcal{N}} ; \boldsymbol{\mu}\right)+\sum_{K \in \mathcal{T}_{h}} \delta_{K} \int_{K}\left(4 x_1(1-x_1) \partial_{x_0} p^{\mathcal{N}}\right)\left(h_{K} \partial_{x_0} z^{\mathcal{N}}\right), \quad z^{\mathcal{N}}, p^{\mathcal{N}} \in Y^{\mathcal{N}}.
   \end{aligned}
\end{equation*}
We consider a parameter space $\mathcal{P}:=\big[10^4,10^6\big]$ and a quite coarse mesh of size $h=0.029$ for the FEM spaces. 
The training set $\mathcal{P}_{\text{train}}$ has cardinality $N_{\text{train}}=100$. We choose $\delta_K =1.0$ \emph{for all} $K \in \mathcal{T}_{h}$ and the penalization term is $\alpha=0.01$. We pursue the convergence in the $L^2$-norm of the state to the desired solution profile $y_d(x)=1.0$, function defined on $\Omega_{obs}$ of Figure \ref{fig:Geometry-Graetz}.
We perform the POD algorithm for $N_{\max}=20$ in a partitioned approach. 
We illustrate the reduced solution for the state and adjoint variables in the best {relative error scenario {in Figures \ref{fig:onoff-graetz-yp1} and \ref{fig:onoff-graetz-yp2}}. Namely, we plot the Only-Offline and Online-Offline Stabilized solutions for $N=1$ and $N=6$. The values of $N$ can be deduced by Figure \ref{fig:error-graetz-yup}. From the gradient equation \eqref{supg-system}, we expect the distributed control $u$ to be equal to the adjoint $p$ up to the multiplicative constant $\alpha$.}
We consider the relative errors between the FEM and the reduced solution in Figure \ref{fig:error-graetz-yup}. We use a testing set $\mathcal{P}_{\text{test}}$ of $100$ elements in $\mathcal{P}$. As previously cited, at $N=6$ we reach the minima for all the three errors for the Online-Offline stabilization; more precisely for the state we touch $e_{y, 6} = 9.65 \cdot 10^{-9}$, for the adjoint $e_{p, 6}= 1.98 \cdot 10^{-8}$ and the control $e_{u, 6}=6.00 \cdot 10^{-9}$. In contrast with this situation, Only-Offline stabilization never falls under $10^{-2}$. This implies that the best choice is to pursue the Online-Offline stabilization procedure for this problem.
However, after $N=6$ the errors begin slightly to increase. Our interpretation to this fact relies on $\mathcal{P}$. Despite the fact that this parameter space might be too large, however the coefficient which multiplies the diffusion operator is still absolutely low in value for every $\mu_1 \in \mathcal{P}$, nearly $10^{-4}$ and $10^{-5}$. Therefore, {also thanks to SUPG stabilization and the distributed control action, the majority of snapshots can be very similar referring to the solution for $\mu_1=10^5$: this translates in  very few  bases to reach a good relative error}. As a matter of fact, the eigenvalues $\lambda^{y}_7, \lambda^{u}_7$ are  $\approx 10^{-15}$ and $\lambda^{p}_7 \approx 10^{-16}$; by their decreasing order, all the subsequent eigenvalues are very close to machine zero.
Thus, recalling \eqref{eigen-const} it follows that all basis components with $N \geq7$ are affected by some rounding errors due to the orthonormalization procedure of the POD (for details see \cite{quarteroni2015reduced}).

Finally, we take a look at the speedup-index in Table \ref{speed-graetz}. All the average values are better for the Only-Offline stabilized ROM procedure due to the fact that the stabilized forms are not taken into account in the online phase. However, the Online-Offline stabilized reduced solution shows very good behaviour, for instance, we have an average equal to $284.3$ for $N=6$. Generally, in this case speedup-index takes average value around $2 \cdot 10^2$ order of magnitude for the first $20$ basis elements.
\begin{table}[!h]
\centering
\begin{tabular}{|c|c|c|c|c|c|c|}
\hline & \multicolumn{3}{c}{Only-Offline Stabilization}& \multicolumn{3}{|c|}{Offline-Online Stabilization} \\
\hline$N$ & \ $\min$ \ & \ average \ & $\max$  & \ $\min$ \ & \ average \ &  $\max$  \\
\hline $1$ & $162.1$ & $296.6$ & $338.1$ & $170.8$ & $261.7$ & $285.9$  \\
$2$  &  $172.2$ & $342.1$ & $391.3$ & $178.4$ & $298.5$ & $327.1$ \\
$3$  &  $168.5$ & $336.2$ & $383.7$  & $192.0$ & $298.9$ & $325.3$\\
$4$  & $165.1$ & $336.1$ & $385.6$ &   $256,7$ & $298.0$ & $322.6$  \\
$5$  & $164.8$ & $331.6$ & $376.6$ &   $220.1$ & $287.6$ & $307.7$  \\
$6$  & $198.3$ & $321.0$ & $366.4$ &   $192.3$ & $284.3$ & $305.7$  \\
$7$  & $186.1$ & $318.4$ & $348.6$ &   $228.6$ & $282.6$ & $306.9$  \\
\hline
\end{tabular}
     \caption{Speedup-index of the Graetz-Poiseuille Problem for Online-Offline and Only-Offline stabilizations with $\mathcal{P}_{\text{test}}$ sampled from $\mathcal{P}=\big[10^4,10^6]$, $N_{\text{test}}=100$, $h=0.029$, $\delta_K =1.0$, $\alpha=0.01$.}
     \label{speed-graetz}
\end{table}
Let us give a look to the {parabolic} version of Problem \eqref{graetz-system} with null initial condition. The {parabolic} Graetz-Poiseuille problem without control has been presented in \cite{pacciarini2014stabilized, torlo2018stabilized}, instead the OCP Graetz Problem under boundary control without Advection-dominancy is studied in \cite{strazzullo2020pod,strazzullo2021model}.

Recalling Figure \ref{fig:Geometry-Graetz}, for a fixed $T>0$ we state the parabolic Graetz-Poiseuille Problem as follows:
\begin{equation} \label{par-graetz-problem}
\begin{cases}
\displaystyle
      \frac{\partial y(\boldsymbol{\mu})}{\partial t}-\frac{1}{{\mu_{1}}} \Delta y(\boldsymbol{\mu})+4 x_1(1-x_1) \partial_{x_0} y(\boldsymbol{\mu})=u, & \text { in } \Omega \times (0,T),\\
\displaystyle
y(\boldsymbol{\mu})=0, & \text { on } \Gamma_{1} \cup \Gamma_{5} \cup \Gamma_{6}  \times (0,T), \\
\displaystyle
y(\boldsymbol{\mu})=1, & \text { on } \Gamma_{2} \cup \Gamma_{4} \times (0,T), \\
\displaystyle
\frac{\partial y(\boldsymbol{\mu})}{\partial \nu}=0, & \text { on } \Gamma_{3} \times (0,T), \\
y(\boldsymbol{\mu})(0)=0, & \text { in } \Omega.
    \end{cases}
\end{equation}
We do simulations in a space-time framework as discussed in Section \ref{time-dep-discr} for a prearranged number of time-steps $N_t$ using a $\mathbb{P}^{1}$-FEM approximation for the high-fidelity solutions. The relative stabilized forms in \eqref{stab-par-block} for derivatives along time for state and adjoint are, respectively:
\begin{equation}
\begin{aligned}
    m_{s}\left(y^{\mathcal{N}}, q^{\mathcal{N}}; \boldsymbol{\mu}\right)&=\left(y^{\mathcal{N}}, q^{\mathcal{N}}\right)_{L^{2}(\Omega)}+\sum_{K \in \mathcal{T}_{h}} \delta_{K} h_{K} \left(y^{\mathcal{N}},   \partial_{x_0} q^{\mathcal{N}}\right)_{K}, \quad y^{\mathcal{N}}, q^{\mathcal{N}} \in Y^{\mathcal{N}},\\
    m^{*}_{s}\left(p^{\mathcal{N}}, z^{\mathcal{N}}; \boldsymbol{\mu}\right)&=\left(p^{\mathcal{N}}, z^{\mathcal{N}}\right)_{L^{2}(\Omega)}-\sum_{K \in \mathcal{T}_{h}} \delta_{K} h_{K}\left(p^{\mathcal{N}}, \partial_{x_0}  z^{\mathcal{N}}\right)_{K}, \quad p^{\mathcal{N}}, z^{\mathcal{N}} \in Y^{\mathcal{N}}.
\end{aligned}
\end{equation}
We consider a final time of $T=3.0$ and a time step of $\Delta t=0.1$, hence we have $N_t = 30$. We choose a quite coarse mesh of $h=0.038$ and the overall high-fidelity dimension is $\mathcal{N}_{tot}=314820$. This means that a single FEM space for a fixed $t$ has a dimension of $\mathcal{N}=3498$. We consider a initial condition of $y_0(x)=0$ \emph{for all} $x \in \Omega$ referring to Figure \ref{fig:Geometry-Graetz}. We want the state solution to converge in the $L^2$-norm to a desired solution profile $y_d(x,t)=1.0$, function defined \emph{for all} $t \in [0,3.0]$ and \emph{for all} $x$ in $\Omega_{obs}$. Here the SUPG stabilization is implemented with parameters $\delta_K =1.0$ \emph{for all} $K \in \mathcal{T}_{h}$. 
$\mathcal{P} := \big[10^4,10^6\big]$ and we choose a training set $\mathcal{P}_{\text{train}}$ of cardinality $N_{\text{train}}=100$. Then, we performed the POD algorithm with $N_{\max}=20$. The penalization parameter is $\alpha=0.01$. 
\begin{figure}[!h]
\centering
     \includegraphics[scale=0.166]{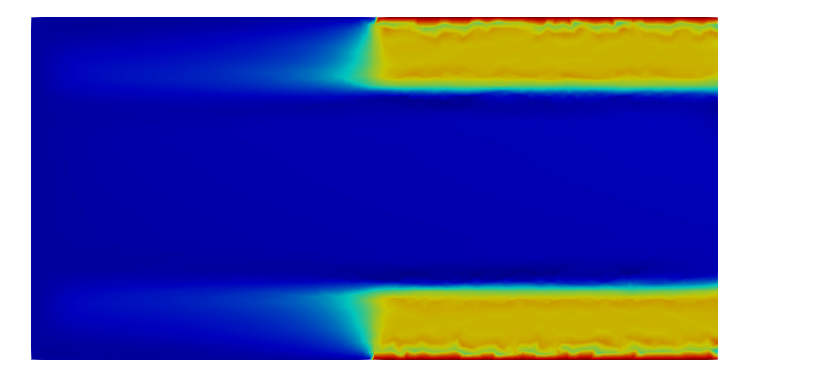}
      \includegraphics[scale=0.166]{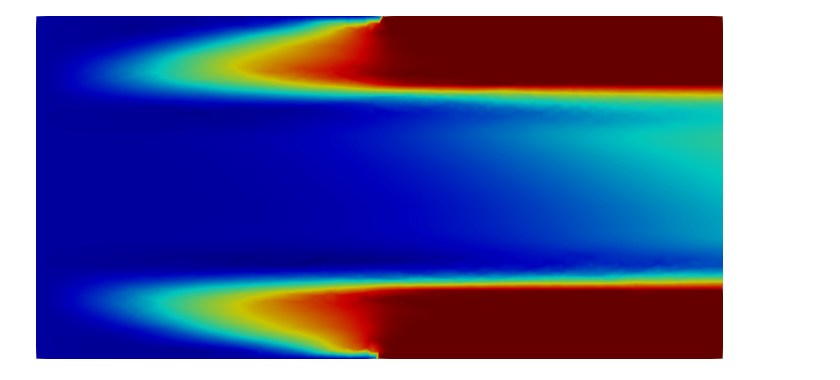}
      \includegraphics[scale=0.166]{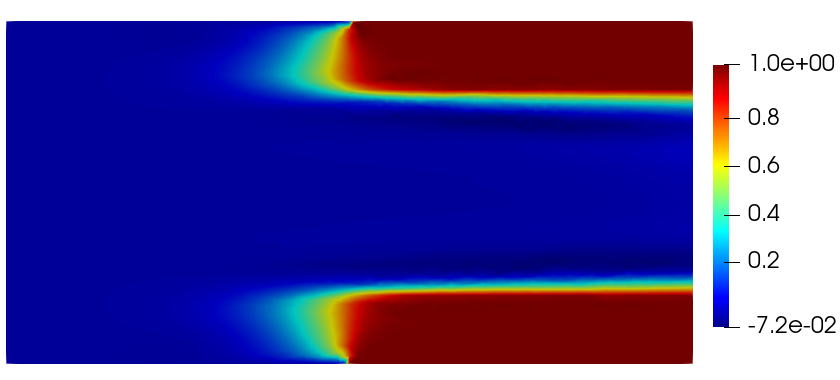} \\
      \includegraphics[scale=0.166]{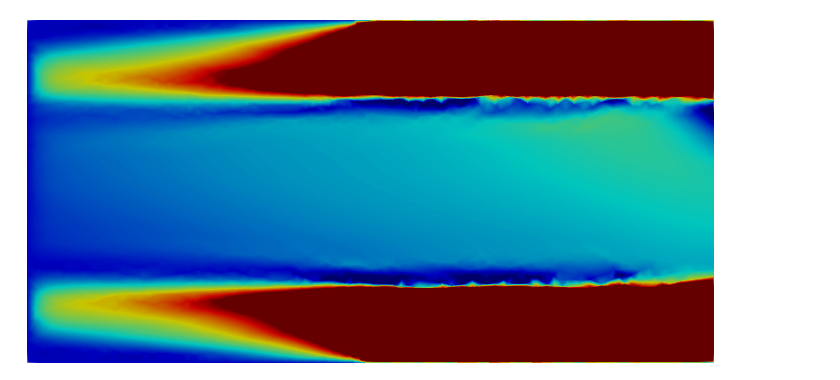}
      \includegraphics[scale=0.166]{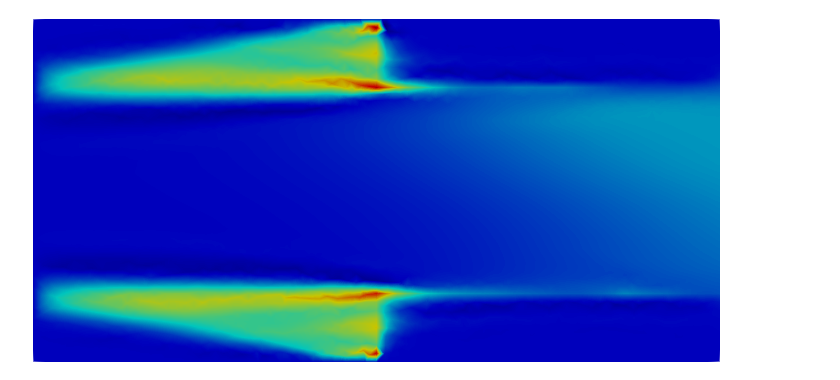}
      \includegraphics[scale=0.166]{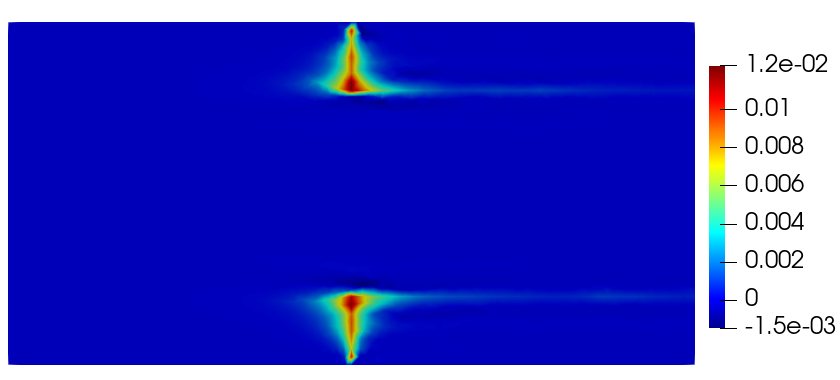} 
     \caption{(\underline{Top}) SUPG FEM solution for the state and (\underline{Bottom}) for the adjoint at $t=0.1$, $t=1.5$, $t=3.0$. {Parabolic} Graetz-Poiseuille Problem, $\mu_1=10^{5}$, $N_t=30$, $h=0.038$, $\delta_K =1.0$, $\alpha=0.01$.}
     \label{fig:par-graetz-fem}
\end{figure}

\begin{figure}[!h]
\centering
     \includegraphics[scale=0.127]{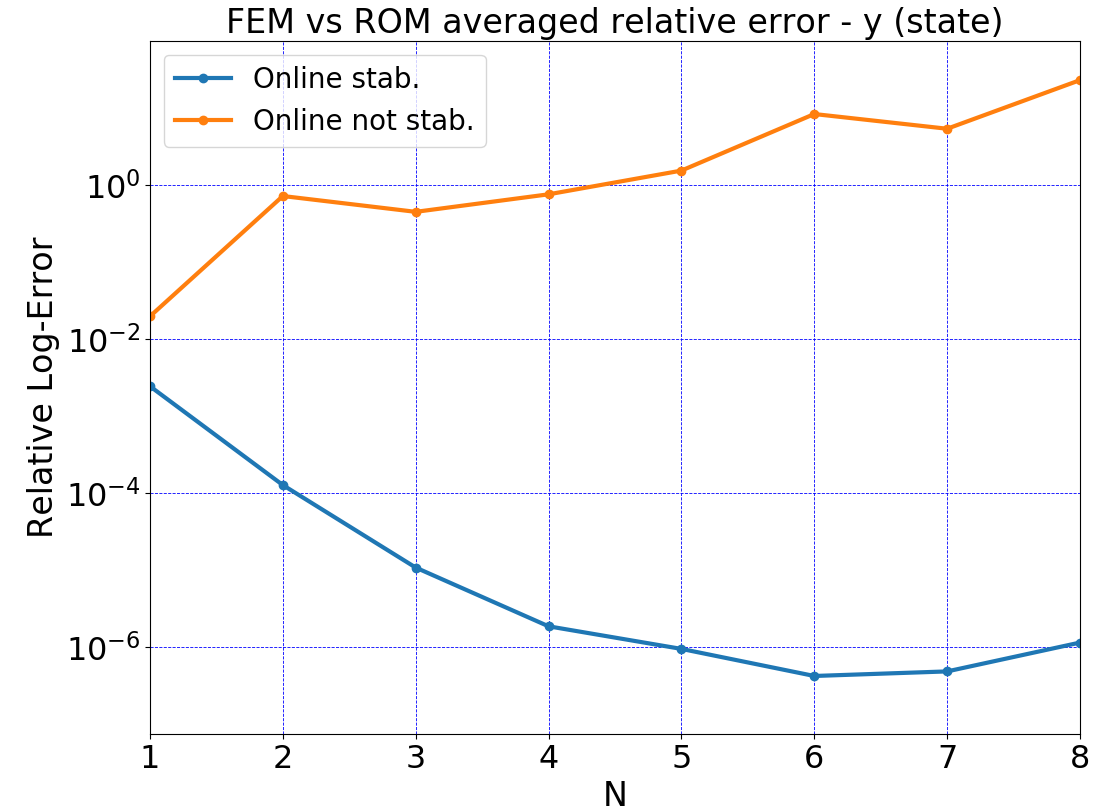} 
      \includegraphics[scale=0.127]{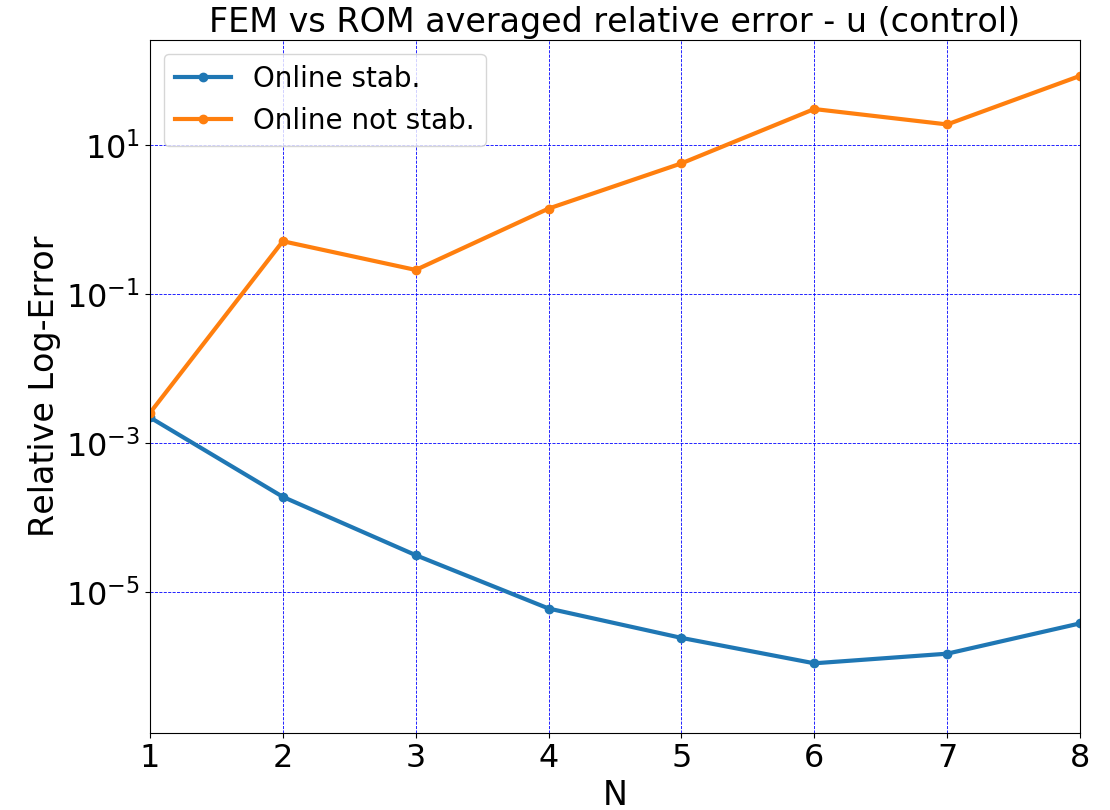}
      \includegraphics[scale=0.127]{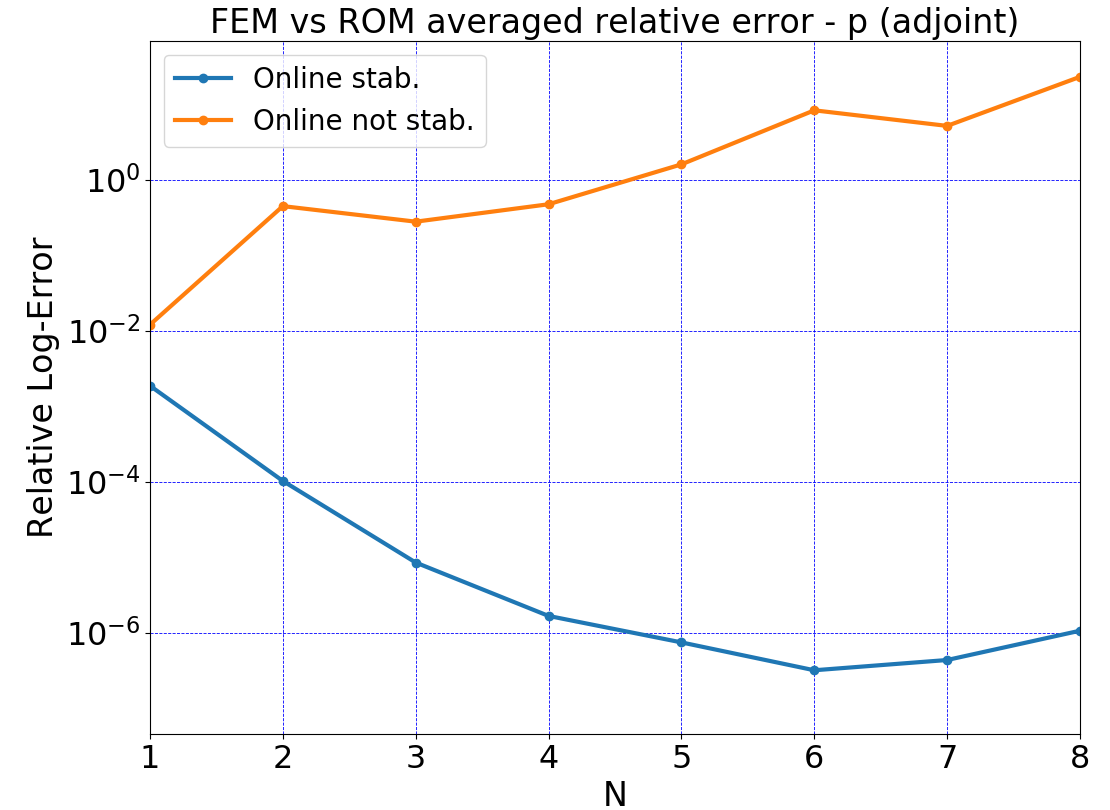}
     \caption{Relative errors between the FEM and Only-Offline and Online-Offline stabilized solutions for the state (left), control (center) and adjoint (right), {Parabolic} Graetz-Poiseuille problem, $N_t=30$, $N_{\text{test}}=100$, $\mathcal{P} = \big[10^4,10^6\big]$, $h=0.038$.}
     \label{fig:error-par-graetz-yup}
\end{figure}

As we will see in Figure \ref{fig:error-par-graetz-yup} the performance of the Only-Offline stabilized reduced solutions are not so good in terms of accuracy, unlike the Online-Offline stabilized ones. 
We consider a testing set of $100$ elements in $\mathcal{P}$.
As succeeded in the steady case in Section \ref{sec:c graetz}, after nearly $N=6$ Online-Offline stabilized errors begin to fluctuate due to the nature of the eigenvalues of the correlation matrix \eqref{correlation-matrix} that are closed to zero machine. For this reason we present the trend of error from $1$ to $10$. However, errors stays close to $10^{-7}$ for the state and the adjoint and $10^{-6}$ to the control. For $N=6$ we have $e_{y, 6}=4.20 \cdot 10^{-7}$, $e_{u, 6}=1.10 \cdot 10^{-6}$ and $e_{p, 6}=3.18 \cdot 10^{-7}$, instead for $N=20$ we have $e_{y, 20}=1.93 \cdot 10^{-7}$, $e_{u, 20}=3.25 \cdot 10^{-7}$ and $e_{p, 20}=1.21 \cdot 10^{-7}$ for the Online-Offline stabilization ROM. 

Finally, the speedup-index trend is shown in Table \ref{speed-graetz-p}. We manage to compute a huge number of reduced solutions in the time of a high-fidelity one: for the Offline-Online stabilization we have an average of nearly $26000$ for $N=6$. On the whole, average speedup-index has an order of magnitude of $2 \cdot 10^4$ for $N \leq 20$.
\begin{table}[!h] 
\centering
\begin{tabular}{|c|c|c|c|c|c|c|}
\hline & \multicolumn{3}{c}{Only-Offline Stabilization}& \multicolumn{3}{|c|}{Offline-Online Stabilization} \\
\hline$N$ & \ $\min$ \ & \ average \ & $\max$  & \ $\min$ \ & \ average \ &  $\max$  \\
\hline $1$ & $21588.3$ & $26588.8$ & $30971.5$ & $18968.4$ & $23588.0$ & $27062.7$  \\
\hline $2$ & $23821.3$ & $29723.4$ & $34817.2$ & $20757.2$ & $26018.9$ & $29929.1$  \\
\hline $3$ & $23571.0$ & $29468.6$ & $34349.5$ & $20547.9$ & $25698.2$ & $29662.5$ \\
\hline $4$ & $23062.2$ & $28880.6$ & $33702.7$ & $21385.2$ & $25380.9$ & $28883.3$  \\
\hline $5$ & $25762.9$ & $28767.9$ & $33488.8$ & $23021.5$ & $25882.4$ & $29388.9$ \\
\hline $6$ & $27003.2$ & $29707.7$ & $34544.7$ & $23236.5$ & $26054.7$ & $29677.5$  \\
\hline $7$ & $26658.5$ & $29481.1$ & $34277.3$ & $23206.6$ & $25879.5$ & $29505.4$ \\
\hline
\end{tabular}
      \caption{Speedup-index of the {Parabolic} Graetz Problem for Online-Offline and Only-Offline stabilization with $\mathcal{P}=\big[10^4,10^6]$, $\alpha=0.01$, $N_t=30$, $N_{\text{test}}=100$, $h=0.038$.}
      \label{speed-graetz-p}
\end{table}

\subsection{Numerical Experiments for Propagating Front in a Square Problem} \label{sec:c square}
\begin{figure}[h!]
\vspace{-.3cm}
   \centering
        \begin{tikzpicture}[scale=3.0]
                  \draw[red] (0,0) -- (0,0.25) node[midway, left, scale=1.2]{$\Gamma_1$};
                 \draw[red] (1,0) -- (0,0) node[midway, below, scale=1.2]{$\Gamma_2$};
                 \draw[blue] (1,1) -- (1,0) node[midway, right, scale=1.2]{$\Gamma_3$};
                 \draw[blue] (0,1) -- (1,1) node[midway, xshift=0.5cm, above, scale=1.2]{$\Gamma_4$};
                 \draw[blue] (0,0.25) -- (0,1) node[midway, left, scale=1.2]{$\Gamma_5$};
                 \draw[black] (0.5,0.45) node[scale=1.4]{$\Omega$};
		         \draw[color=DEblue!100, fill=DEblue!10] (.25, 1) -- (0.25, 0.75) -- (1,0.75) -- (1,1) -- (.25, 1)node[midway, left, scale=1.2]{};
                 \draw[black] (0.65,0.875) node[scale=1.]{{$\Omega_{obs}$}};
                 \draw[black] (0.5,0.5) node[scale=1.4]{};
                 \filldraw[black] (0,0.25) circle (0.3pt) node[above, left]{(0,0.25)};
                 \filldraw[black] (0,1) circle (0.3pt) node[left]{(0,1)};
                 \filldraw[black] (1,0.75) circle (0.3pt) node[right]{(1,0.75)};
                 \filldraw[black] (1,1) circle (0.3pt) node[right]{(1,1)}; \filldraw[black] (0.25,1) circle (0.3pt) node[above]{(0.25,1)};
                 \filldraw[black] (1,0) circle (0.3pt) node[below, right]{(1,0)};
                 \filldraw[black] (0,0) circle (0.3pt) node[below]{(0,0)};
        \end{tikzpicture} 
        \caption{Geometry of the Square Problem}
        \label{fig:geometry-square}
\end{figure}
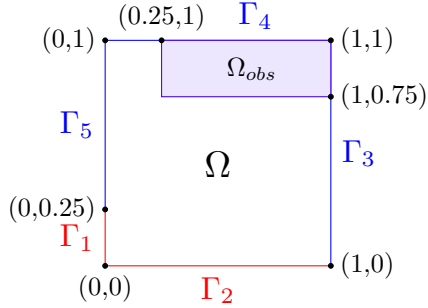

In this Section, we consider a problem studied in the Advection-Dominated form in \cite{pacciarini2014stabilized,torlo2018stabilized} from a numerical point of view and we will add a distributed control to it.
Let $\Omega$ be the unit square in $\mathbb{R}^{2}$. We consider the representation in Figure \ref{fig:geometry-square}.
Also in this case, $(x_0,x_1)$ are the coordinates of the square domain. Referring to Figure \ref{fig:geometry-square}, $\Gamma_1 := \{0\} \times [0,0.25]$, $\Gamma_2 := [0,1] \times \{0\}$, $\Gamma_3 := \{1\} \times [0,1]$, $\Gamma_4 := [0,1] \times \{1\}$, $\Gamma_5 := \{0\} \times [0.25,1]$; $\Omega_{obs}:= [0.25,1]\times [0.75,1]$.
Given $\boldsymbol{\mu} = (\mu_1, \mu_2)$, the problem is formulated as
\begin{equation}\label{square-system}
\begin{cases}
\displaystyle
-\frac{1}{\mu_1} \Delta y(\boldsymbol{\mu})+(\cos{\mu_2},\sin{\mu_2}) \cdot \nabla y(\boldsymbol{\mu})=u, & \text { in } \Omega, \\
\displaystyle
y(\boldsymbol{\mu})=1, & \text { on } \Gamma_{1} \cup \Gamma_{2}, \\
\displaystyle
y(\boldsymbol{\mu})=0, & \text { on } \Gamma_{3} \cup \Gamma_{4} \cup \Gamma_{5}.
\end{cases}
\end{equation}
{We assume{the embedding between $Y=H^{1}(\Omega)$ and $Z=L^2(\Omega_{obs})$} as the Observation operator}.
In our test cases, $\mathcal{P}:=\big[10^4,10^5\big] \times \big[0,1.57\big]$. 
In this case, we have that {the domain of definition of our state $y$ is $$\tilde{Y}:= \big\{v \in H^{1} \big(\Omega\big) \text{ s.t. } \mathrm{BC} \text{ in }  (\ref{square-system}) \big\}.$$ 
Exactly as done in the previous paragraph, we define a lifting function $R_y \in H^1\big(\Omega\big)$ such that satisfies $\mathrm{BC} \text{ in }$  (\ref{square-system}). We define $\bar{y} := y - R_y$, even though we denote $\bar{y}$ as $y$ again for the sake of notation. We consider $Y:= H^1_{0}(\Omega)$, }$U = \mathrm{L}^{2}(\Omega)$ and $Q := Y^{*}$, hence the adjoint $p$ is such that $p=0$ on $\partial \Omega$.
We define the objective functional $J$ exactly as in \eqref{Jgraetz}; instead, $a$ and $b$ are 
\begin{equation*} \label{square-forms}
\displaystyle
a\left(y, p ; \boldsymbol{\mu}\right) := \int\limits_{\Omega} \frac{1}{\mu_{1}} \nabla y \cdot \nabla p +(\cos{\mu_2},\sin{\mu_2}) \cdot \nabla y  p \ \mathrm{dx}, \text{ and } b\left(u, p ; \boldsymbol{\mu}\right) := - \int\limits_{\Omega} u p \ \mathrm{dx}.
\end{equation*}
 and $\langle p, f(\boldsymbol{\mu}) \rangle_{Y^{*}Y} = -a\left(R_y,p; \boldsymbol{\mu}\right).$
Then we follow usual discussions of Sections \ref{sec:problem} and \ref{sec:truth_discretization}.

\begin{figure}[!h]
\centering
     \includegraphics[scale=0.151]{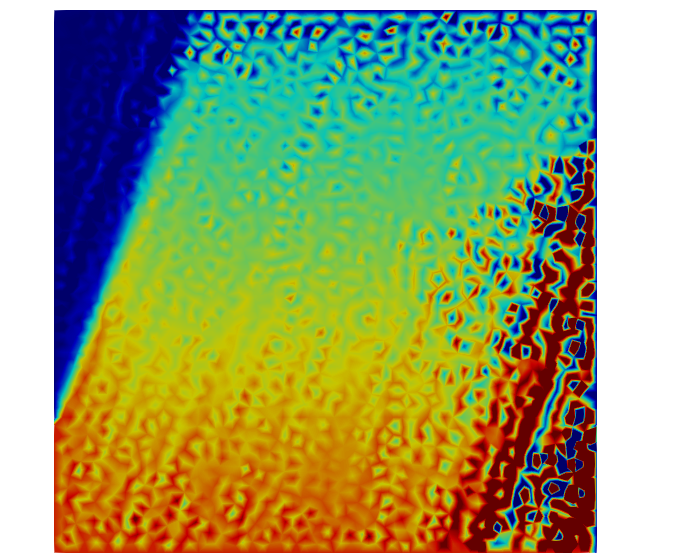} 
     \includegraphics[scale=0.151]{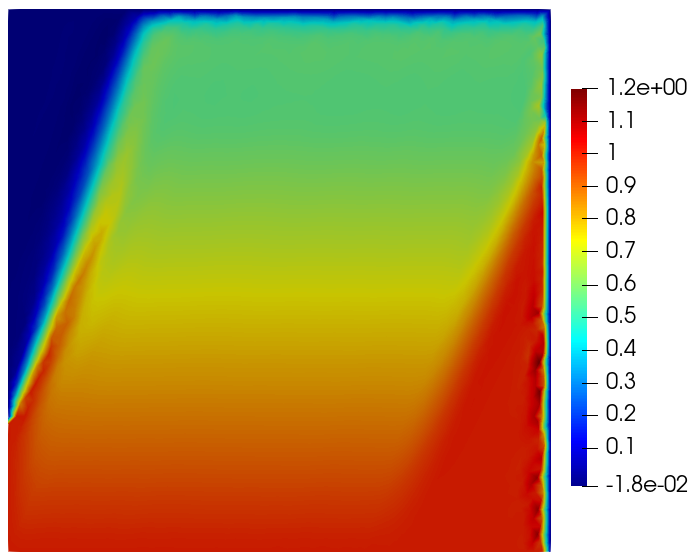}
     \includegraphics[scale=0.151]{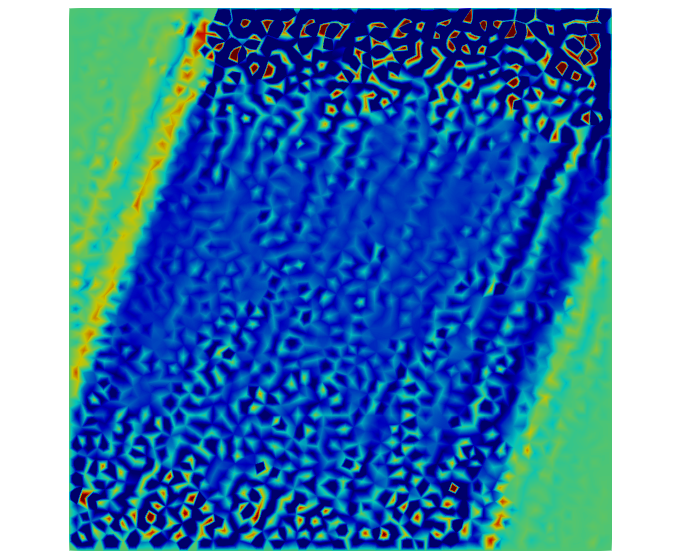} 
     \includegraphics[scale=0.151]{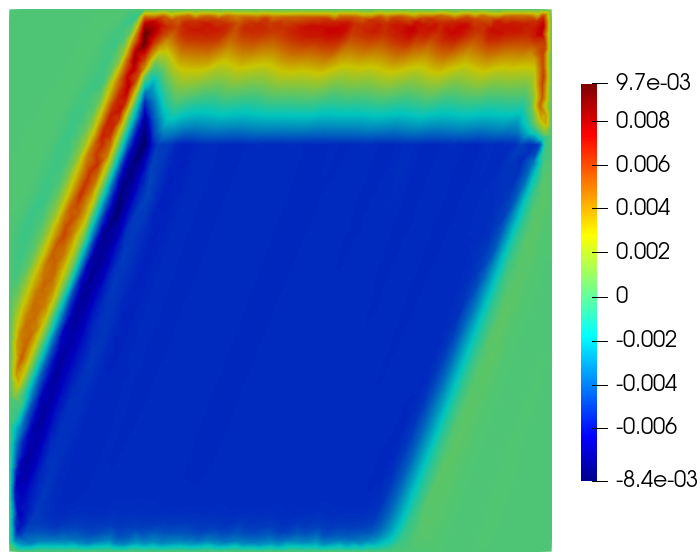}
     \caption{Numerical solution without stabilization and SUPG FEM solution with $\boldsymbol{\mu}=(2 \cdot 10^4,1.2)$ for state (left) and adjoint (right) variables in the Propagating Front in a Square Problem, $\alpha=0.01$, $h=0.025$, $\delta_K =1.0$.}
     \label{fig:fem-square-yup}
\end{figure}

\begin{figure}[!h]
\centering
     \includegraphics[scale=0.151]{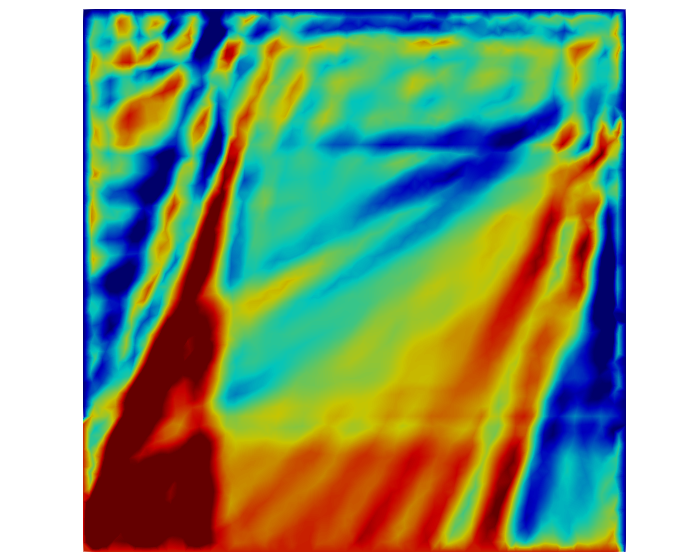}
     \includegraphics[scale=0.151]{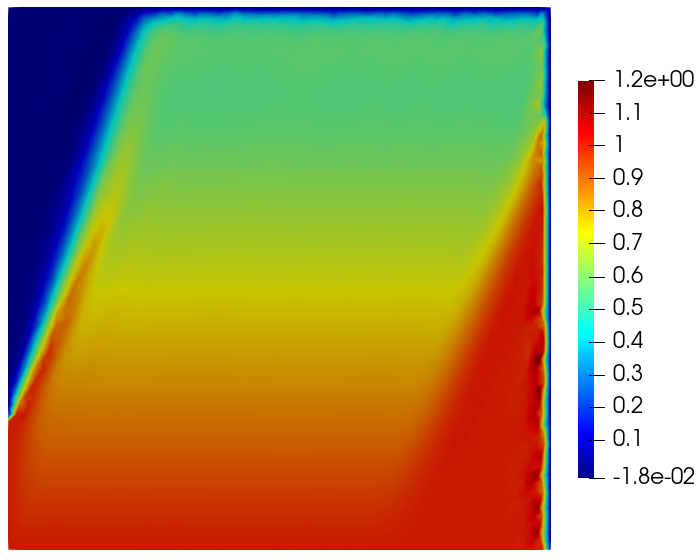}
     \includegraphics[scale=0.151]{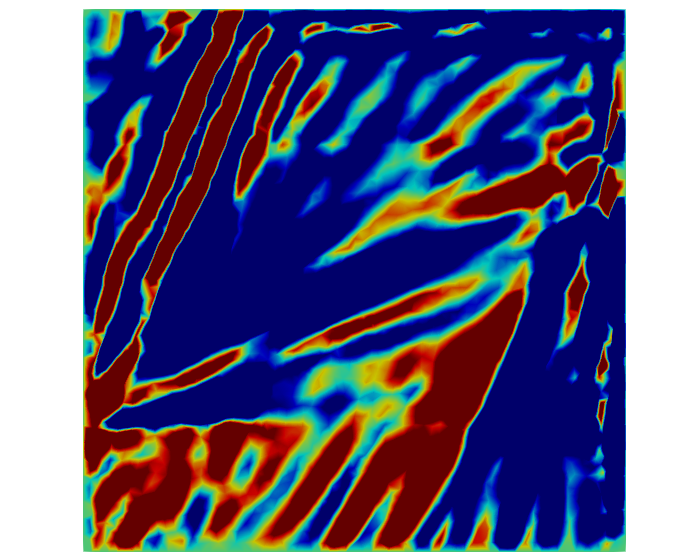}
     \includegraphics[scale=0.151]{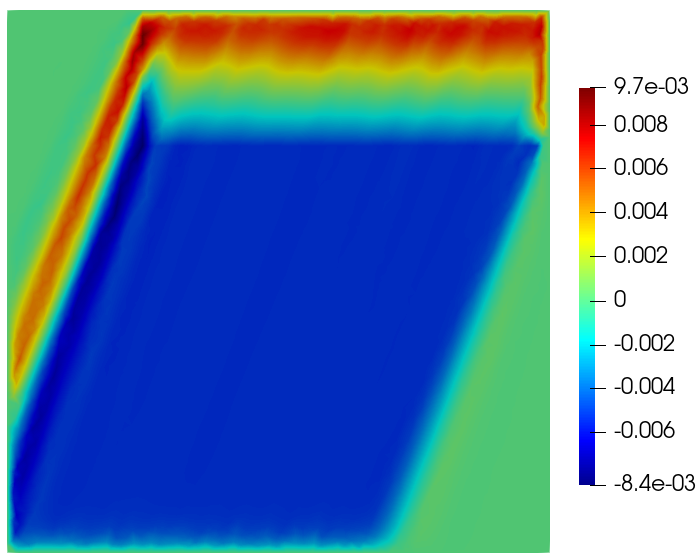}
     \caption{Only-Offline stabilized and Online-Offline stabilized reduced solutions with $\boldsymbol{\mu}=(2 \cdot 10^4,1.2)$ or state (left) and adjoint (right) variables in the Propagating Front in a Square Problem, $\alpha=0.01$, $N=50$, $h=0.025$, $\delta_K =1.0$, $\mathcal{P}=\big[10^4,10^5]\times \big[0,1.57]$.}
     \label{fig:onoff-square-yup}
\end{figure}

\begin{figure}[!h]
\centering
     \includegraphics[scale=0.129]{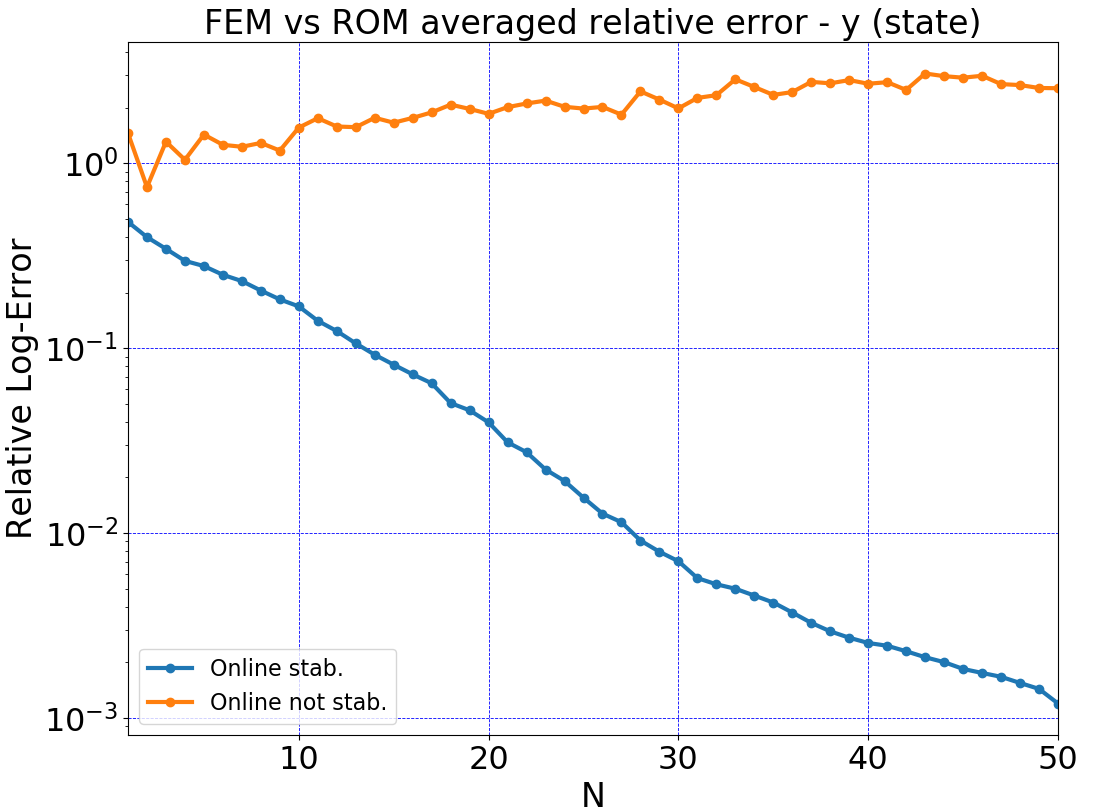} 
     \includegraphics[scale=0.129]{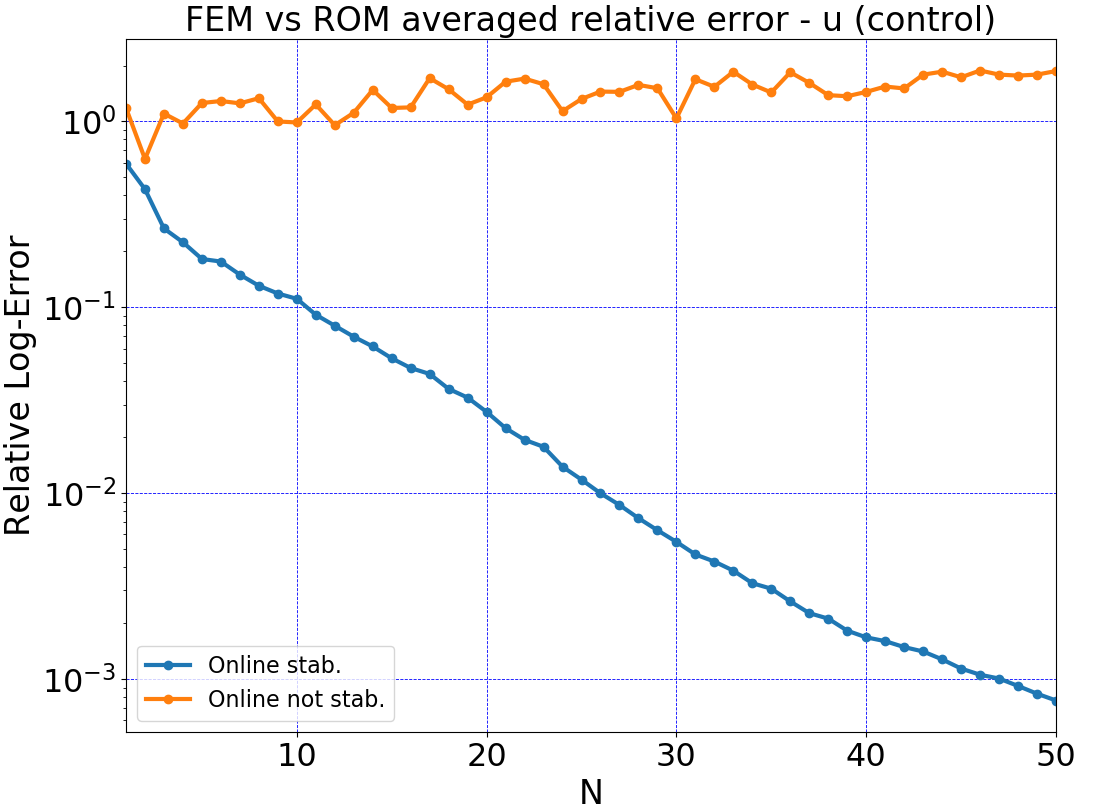}
     \includegraphics[scale=0.129]{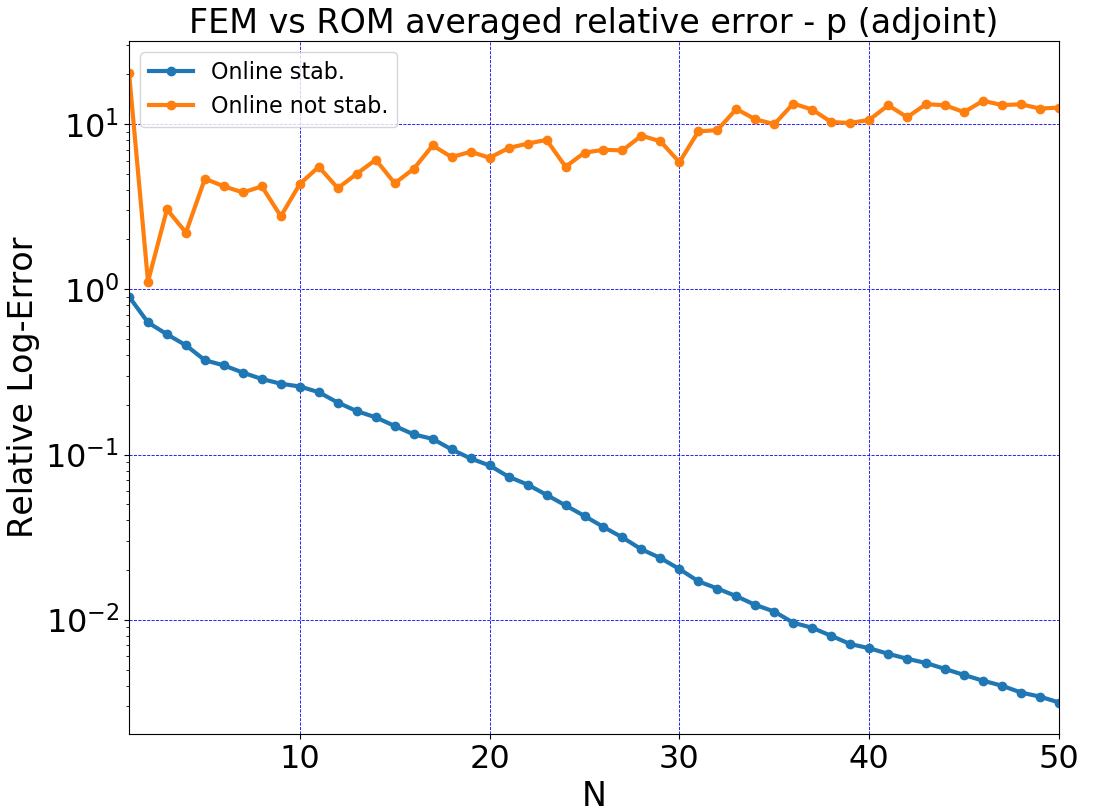}
     \caption{Relative errors between FEM and reduced solutions with $\mathcal{P}=\big[10^4,10^5]\times \big[0,1.57]$ for the state (left), control (center) and adjoint (right) in the Propagating Front in a Square Problem, $N_{\text{test}}=100$, $\alpha=0.01$, $h=0.025$, $\delta_K =1.0$.}
     \label{fig:error-square-yup}
\end{figure}

We exploit a $\mathbb{P}^1$-FEM approximation for the optimality system by using the usual SUPG stabilization technique, arriving to system \eqref{supg-system}. We remark that the stabilized forms $a_s$  and $a^{*}_s$ are, respectively:
\begin{equation*}
  \begin{aligned}
a_{s}\left(y^{\mathcal{N}}, q^{\mathcal{N}} ; \boldsymbol{\mu}\right)&:=a\left(y^{\mathcal{N}}, q^{\mathcal{N}} ; \boldsymbol{\mu}\right)+\sum_{K \in \mathcal{T}_{h}} \delta_{K} \int_{K}h_{K} \left( \cos \mu_2, \sin \mu_2 \right) \cdot  \nabla y^{\mathcal{N}} \left( \cos \mu_2, \sin \mu_2 \right) \cdot  \nabla q^{\mathcal{N}}, \\
   a^{*}_{s}\left(z^{\mathcal{N}}, p^{\mathcal{N}} ; \boldsymbol{\mu}\right)&:=a^{*}\left(z^{\mathcal{N}}, p^{\mathcal{N}} ; \boldsymbol{\mu}\right)+\sum_{K \in \mathcal{T}_{h}} \delta_{K} \int_{K}h_{K} \left( \cos \mu_2, \sin \mu_2 \right) \cdot  \nabla p^{\mathcal{N}} \left( \cos \mu_2, \sin \mu_2 \right) \cdot  \nabla z^{\mathcal{N}},
   \end{aligned}
\end{equation*}
\emph{for all} $y^{\mathcal{N}}, q^{\mathcal{N}}, z^{\mathcal{N}}, p^{\mathcal{N}} \in Y^{\mathcal{N}}$. As previously done, we build a training set $\mathcal{P}_{\text{train}}$ and a testing set $\mathcal{P}_{\text{test}}$ with both cardinality $n_{\text{train}}=100$. The mesh size $h$ is $0.025$ and therefore the overall dimension of the high-fidelity approximation is $12087$, which implies that state, control and adjoint spaces have dimension equal $\mathcal{N}=4029$. The SUPG stabilization is implemented with parameters $\delta_K =1.0$ \emph{for all} $K \in \mathcal{T}_{h}$. The penalization parameter is $\alpha=0.01$ and we pursue {the state solution to be convergent in the $L^2$-norm to a} desired solution profile $y_d(x)=0.5$, defined \emph{for all} $x$ in $\Omega_{obs}$ of Figure \ref{fig:geometry-square}. 
In Figure \ref{fig:fem-square-yup} we observe state and adjoint FEM solutions for $\boldsymbol{\mu}=(2\cdot 10^4, 1.2)$. Instead, in Figure \ref{fig:onoff-square-yup} we illustrate Only-Offline and Online-Offline reduced solution for the state and the adjoint variable with $\boldsymbol{\mu}=(2\cdot 10^4, 1.2)$ for $N=50$. 
\quad
These computational evidences and the analysis of the relative errors {show that Online-Offline stabilization procedure is preferable in this setting}. In Figure \ref{fig:error-square-yup}, the trend is the same of all three variables, where errors continue to decrease along all $N$: we have $e_{y, 50} = 1.20 \cdot 10^{-3}$, $e_{u, 50} = 7.67 \cdot 10^{-4}$ and $e_{p, 50} = 3.16 \cdot 10^{-3}$.
Concerning the speedup-index, the performance are quite good as seen in Table \ref{speed-s}. For the best approximation, we have that we can compute an average of $44$ Online-Offline reduced solutions when we build the associated FEM one. Obviously, the Only-Offline stabilized one is slightly better. On the whole, speedup-index takes average value around $10^1,10^2$ order of magnitude for $N\leq 50$.
\begin{table}[!h] 
\centering
\begin{tabular}{|c|c|c|c|c|c|c|}
\hline & \multicolumn{3}{c}{Only-Offline Stabilization}& \multicolumn{3}{|c|}{Offline-Online Stabilization} \\
\hline$N$ & $\min$ & average & $\max$ & $\min$ & average& $\max$  \\
\hline $1$ & $123.1$ & $198.8$ & $243.1$ & $110.0$ & $162.7$ & $181.9$  \\
$10$  &  $132.2$ & $200.4$ & $244.2$ & $110.8$ & $158.3$ & $176.9$ \\
$20$   &  $84.6$ & $158.3$ & $191.7$  & $60.1$ & $124.3$ & $141.3$\\
$30$  & $78.8$ & $114.7$ & $137.8$ &   $65.0$ & $92.5$ & $104.7$  \\
$40$  & $54.2$ & $78.6$ & $96.1$ &   $46.9$ & $64.2$ & $72.3$  \\
$50$  & $33.2$ & $53.0$ & $64.8$ &   $28.5$ & $44.0$ & $49.9$  \\
\hline
\end{tabular}
     \caption{Speedup-index of the Propagating Front in a Square Problem for Online-Offline and Only-Offline stabilization with training set $\mathcal{P}=\big[10^4,10^5]\times \big[0,1.57]$, $\alpha=0.01$, $N_{\text{test}}=100$, $h=0.025$, $\delta_K =1.0$.}
     \label{speed-s}
\end{table}
Now we study the {parabolic} case of the Propagating Front in a Square Problem for a fix $T>0$:
\begin{equation}\label{par-square-problem}
\begin{cases}
\displaystyle
\frac{\partial y(\boldsymbol{\mu})}{\partial t}-\frac{1}{\mu_1} \Delta y(\boldsymbol{\mu})+(\cos{\mu_2},\sin{\mu_2}) \cdot \nabla y(\boldsymbol{\mu})=u, & \text { in } \Omega  \times (0,T), \\
\displaystyle
y(\boldsymbol{\mu})=1, & \text { on } \Gamma_{1} \cup \Gamma_{2}  \times (0,T), \\
\displaystyle
y(\boldsymbol{\mu})=0, & \text { on } \Gamma_{3} \cup \Gamma_{4} \cup \Gamma_{5}  \times (0,T), \\
\displaystyle
y(\boldsymbol{\mu})(0)=0, & \text { in } \Omega,\\
\end{cases}
\end{equation}
with initial value
$y_0(x)=0$ \emph{for all} $x \in {\Omega}$ referring to the domain in Figure \ref{fig:geometry-square}. We consider a final time $T=3.0$ and a time-step $\Delta t=0.1$, hence $N_t=30$. We choose a quite coarse mesh of size $h=0.036$ and the overall dimension of the space-time system is $N_{tot}=174780$, which means that a single FEM space has dimension $\mathcal{N}=1942$. {Again, our aim is to achieve in a $L^2$-mean a desired solution profile $y_d(x,t)=0.5$, defined \emph{for all} $t \in [0,3]$ and $x$ in $\Omega_{obs}$ of Figure \ref{fig:geometry-square}.  The penalization parameter is $\alpha=0.01$. 
We set $\delta_K =1.0$ \emph{for all} $K \in \mathcal{T}_{h}$. Here the stabilized forms in \eqref{stab-par-block} for state and adjoint equations are, respectively:
\begin{equation*}
\begin{aligned}
    m_{s}\left(y^{\mathcal{N}}, q^{\mathcal{N}}; \boldsymbol{\mu}\right)&=\left(y^{\mathcal{N}}, q^{\mathcal{N}}\right)_{L^{2}(\Omega)}+\sum_{K \in \mathcal{T}_{h}} \delta_{K} h_{K} \left(y^{\mathcal{N}},  \left( \cos \mu_2, \sin \mu_2 \right) \cdot  \nabla q^{\mathcal{N}}\right)_{K}, \quad y^{\mathcal{N}}, q^{\mathcal{N}} \in Y^{\mathcal{N}},\\
    m^{*}_{s}\left(p^{\mathcal{N}}, z^{\mathcal{N}}; \boldsymbol{\mu}\right)&=\left(p^{\mathcal{N}}, z^{\mathcal{N}}\right)_{L^{2}(\Omega)}-\sum_{K \in \mathcal{T}_{h}} \delta_{K} h_{K}\left(p^{\mathcal{N}}, \left( \cos \mu_2, \sin \mu_2 \right) \cdot  \nabla z^{\mathcal{N}}\right)_{K}, \quad p^{\mathcal{N}}, z^{\mathcal{N}} \in Y^{\mathcal{N}}.
\end{aligned}
\end{equation*}
We consider a parameter space equal to the steady case, i.e.\ $\mathcal{P} := \big[10^4,10^5\big] \times \big[0,1.57\big]$. Our training set has cardinality $N_{\text{train}}=100$. In Figure \ref{fig:par-fem-square-y-p} we show a representative stabilized FEM solution for $\boldsymbol{\mu}=(2\cdot 10^{4},1.2)$ for some instants of time. We choose to perform a POD procedure with $N_{max}=30$.
\begin{figure}[!h]
\centering
     \includegraphics[scale=0.18]{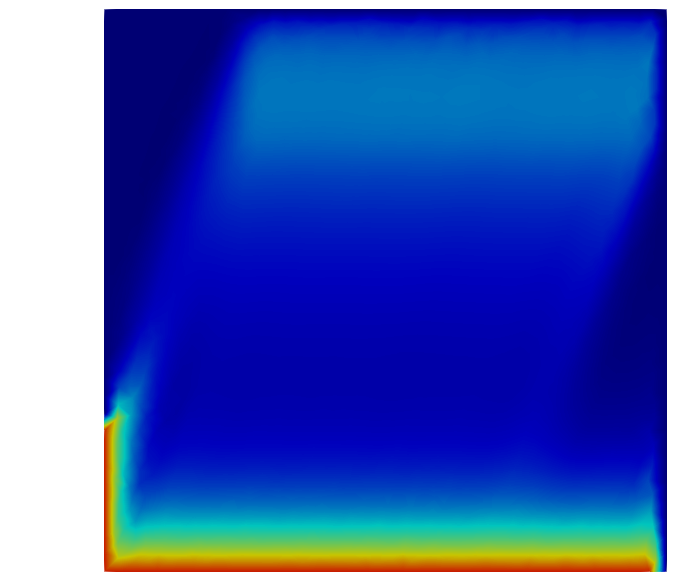}
     \includegraphics[scale=0.18]{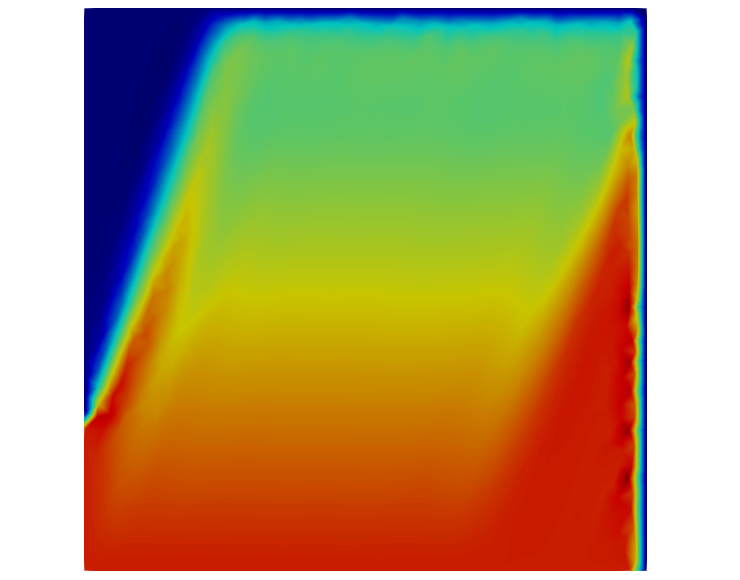}
     \includegraphics[scale=0.18]{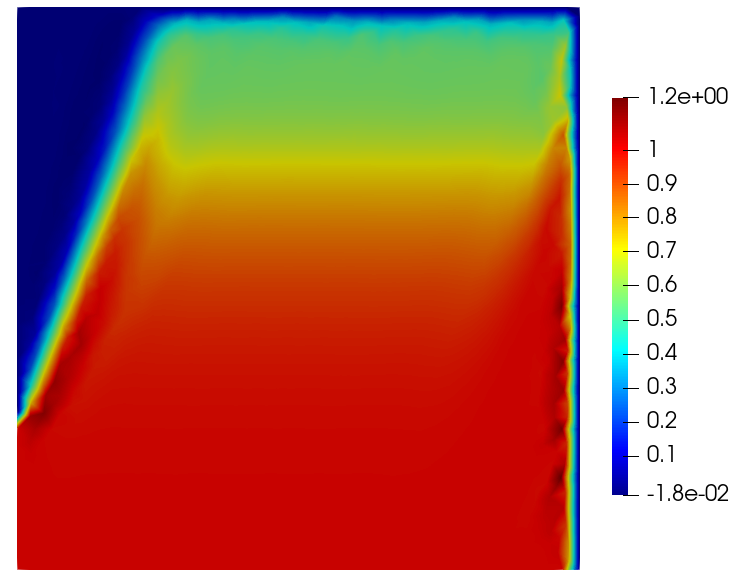}
     \includegraphics[scale=0.18]{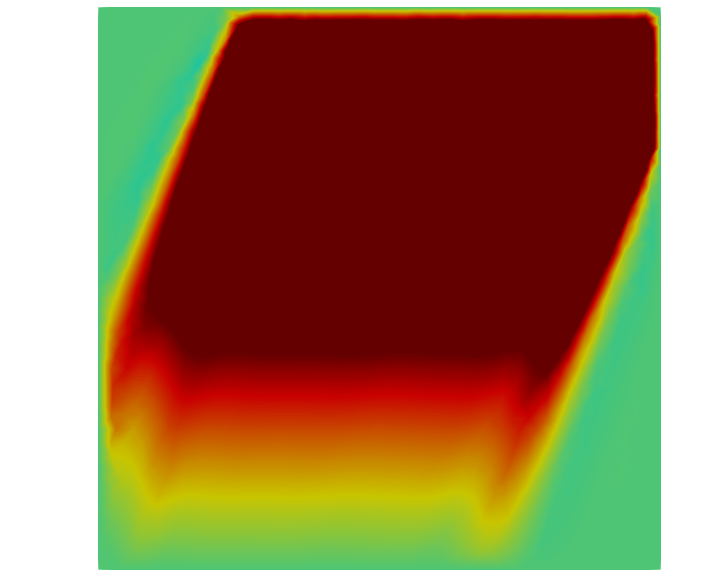}
     \includegraphics[scale=0.18]{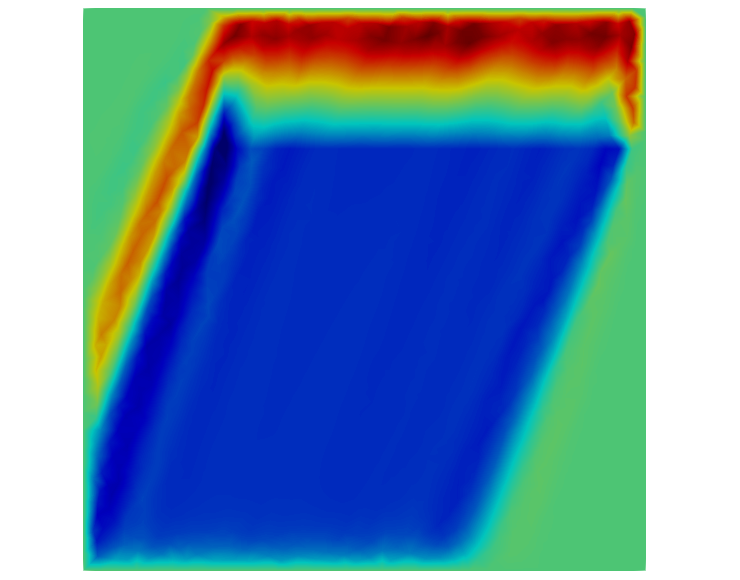}
     \includegraphics[scale=0.18]{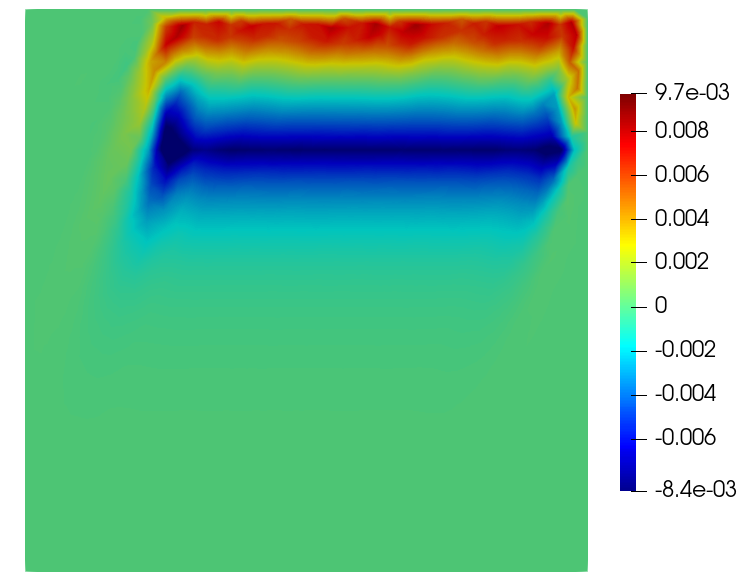}
     \caption{ (\underline{Top}) SUPG FEM state solution and (\underline{Bottom}) SUPG FEM adjoint solution for $\boldsymbol{\mu}=(2\cdot 10^{4},1.2)$ for time $t=0.1$ (left), $t=1.5$ (center), $t=3.0$ (right), in the Parabolic Propagating Front in a Square Problem, $h=0.036$, $\alpha=0.01$, $\delta_K =1.0$.}
     \label{fig:par-fem-square-y-p}
\end{figure}
In Figure \ref{fig:error-par-square-yup} one can see the relative errors of the three variables. As previously said, Only-Offline procedure has not good error behaviour. Instead, it is worth to note that in a Online-Offline stabilization context, errors between the FEM and the reduced solutions decrease as $N$ grows. The fact that we deal with a two-dimensional parameter space implies to require more $N$ basis for a good approximation of the reduced solution. We have $e_{y, 30}=2.17 \cdot 10^{-3}$, $e_{u, 30}=1.59 \cdot 10^{-3}$ and $e_{p, 30}=5.62 \cdot 10^{-3}$. Therefore, also for this case test we can state that the SUPG stabilization is an efficient procedure for the ROMs.
\begin{figure}[!h]
\centering
     \includegraphics[scale=0.13]{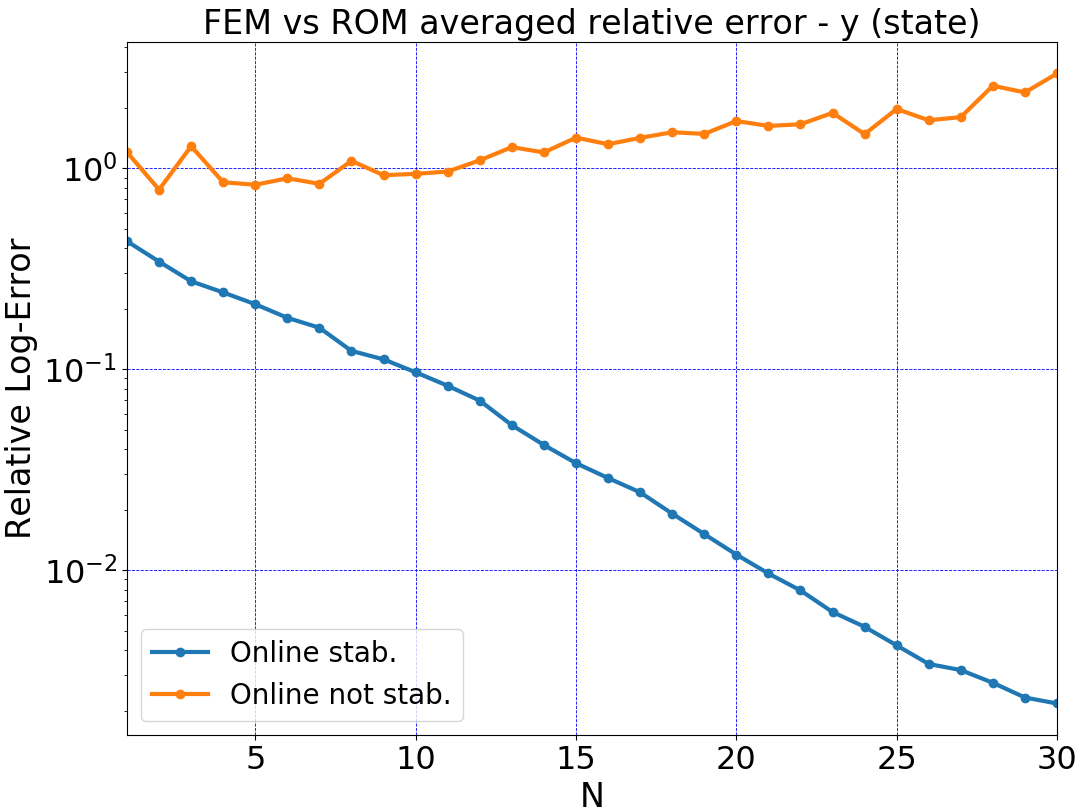} 
      \includegraphics[scale=0.13]{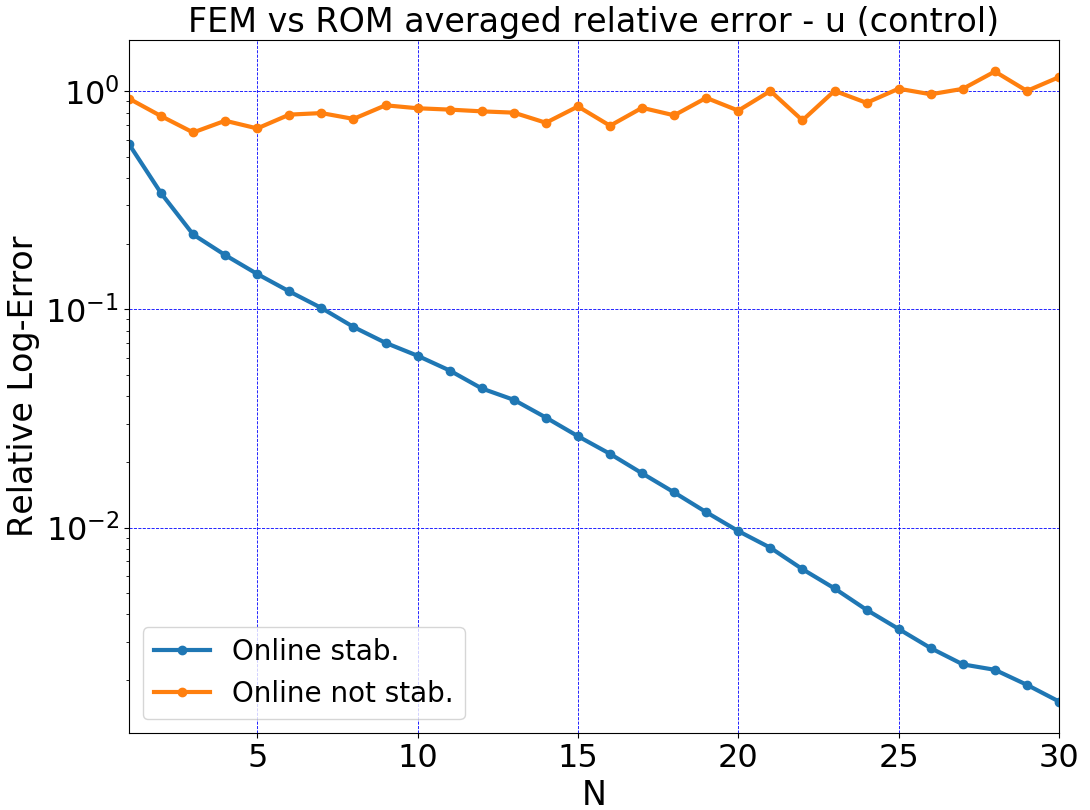} 
      \includegraphics[scale=0.13]{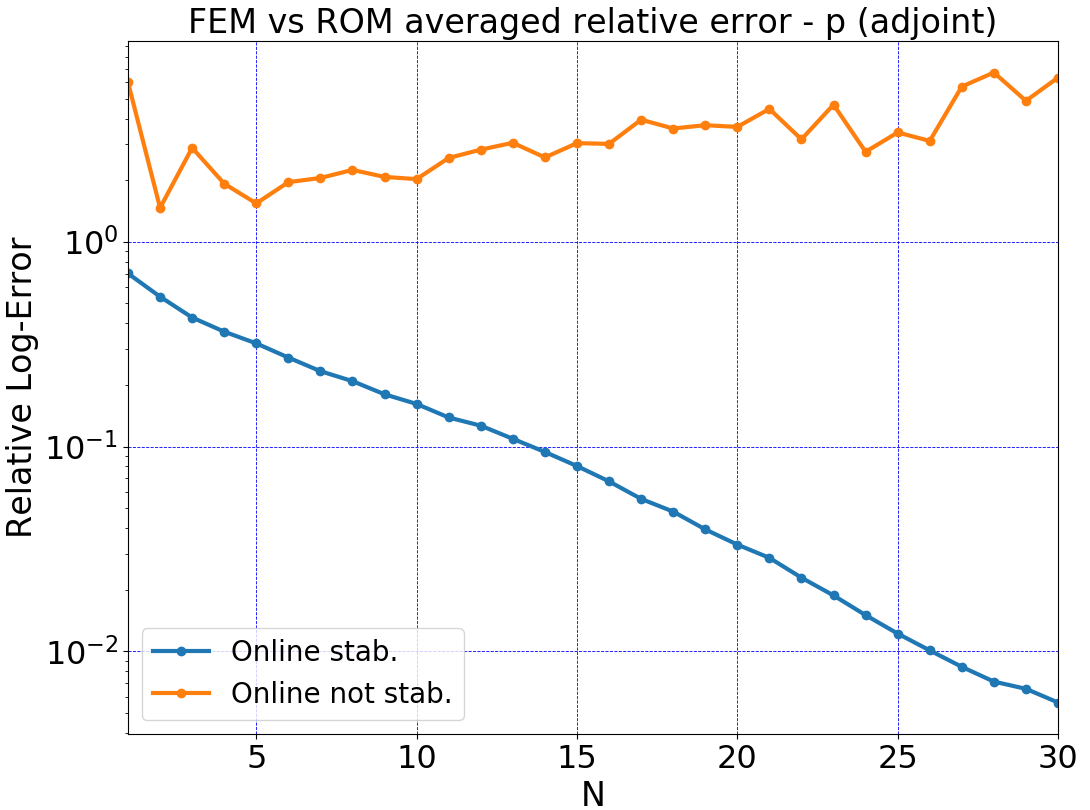} 
     \caption{Relative errors between the FEM and the Only-Offline and Online-Offline stabilized reduced solution for the state (left), control (center) and adjoint (right) solutions, respectively with $\mathcal{P} = \big[10^4,10^5\big] \times \big[0,1.57\big]$, $N_t=30$, $N_{\text{test}}=100$, $\delta_K =1.0$, $h=0.036$.}
     \label{fig:error-par-square-yup}
\end{figure}

Finally, we show the results about the speedup-index in Table \ref{speed-p-s}. For $N=30$, {not only we have the best accuracy for the reduced problem, but we are able to computed, averagely, $3981$ reduced solution in the interval of a FEM simulation}. Speedup-index has an average order of magnitude of $10^3$ overall.
\begin{table}[h!]
\centering
\begin{tabular}{|c|c|c|c|c|c|c|}
\hline & \multicolumn{3}{c}{Only-Offline Stabilization}& \multicolumn{3}{|c|}{Offline-Online Stabilization} \\
\hline$N$ & $\min$ & average & $\max$ & $\min$ & average& $\max$  \\
\hline $5$ & $6136.6$ & $8659.4$ & $12654.5$ & $4588.3$ & $6605.8$ & $9914.7$  \\
$10$  &  $6018.7$ & $8278.8$ & $11989.5$ & $4282.9$ & $6231.7$ & $9353.3$ \\
$15$   &  $5611.3$ & $7721.4$ & $11359.5$  & $3725.3$ & $5779.0$ & $8794.3$\\
$20$  & $4820.0$ & $7041.2$ & $10143.0$ &   $3567.7$ & $5318.9$ & $8091.1$  \\
$25$  & $3814.1$ & $5970.5$ & $9212.2$ &   $2751.7$ & $4366.6$ & $6678.3$  \\
$30$  & $3432.0$ & $5420.1$ & $8462.8$ &   $2424.7$ & $3981.3$ & $6147.9$  \\
\hline
\end{tabular}
  \caption{Speedup-index of the {parabolic} Propagating Front in a Square Problem for Online-Offline and Only-Offline stabilization with training set $\mathcal{P}:=\big[10^4,10^5]\times \big[0,1.57]$, $h=0.036$, $\alpha=0.01$, $N_t=30$, $N_{\text{test}}=100$, $\delta_K =1.0$.}
  \label{speed-p-s}
\end{table}

\section{Conclusions and Perspectives}
In this work, we presented the numerical experiments concerning Advection-Dominated OCPs in a ROM context with high P\'eclet number, both in the steady and the {parabolic} cases, under SUPG stabilization. Concerning ROMs, we can have two possibilities of stabilization: we can apply SUPG only to the offline phase or we can use it in both online and offline phases. We analyzed relative errors between the reduced and the high fidelity solutions and of the speedup-index concerning the Graetz-Poiseuille and Propagating Front in a Square Problems, always under a distributed control. 

A $\mathbb{P}^{1}$-FEM approximation for the state, control and adjoint spaces is used in a \textit{optimize-then-discretize} framework. For parabolic problems, a space-time approach is followed and we applied in a suitable way the SUPG stabilization. For the ROM, we considered a \textit{partitioned approach} for all three variables using the POD algorithm. In all the steady and unsteady experiments, the ROM technique performed excellently in a Online-Offline stabilization framework. Especially for parabolic problems, the speedup-index features large values thanks to the space-time formulation. Only-Offline stabilization technique performed very poorly in terms of errors, despite the little favorable speedup values.{When the speedup indexes are comparable, the accuracy results determine what strategy to use.} Thus, Online-Offline stabilization is preferable.

We also performed experiments inherent a geometrical parametrization and boundary control for the Graetz-Poiseuille Problem that are not shown here. Results were quite good for Online-Offline stabilization: we had some little oscillations regarding relative errors due to the complexity of the problem. As a perspective, it might be interesting to create a strongly-consistent stabilization technique that flattens all the fluctuation for these two configurations, since, to the best of our knowledge, this topic is still a novelty in literature.

Regarding the SUPG stabilization for parabolic OCPs in a \textit{optimize-then-discretize} framework, it would be also worth deriving some theoretical results that give us the accuracy of the numerical solution with respect to the time-step and the mesh-size.\\ {Moreover, another possibility for future developments would be the analysis of the ``discretize-then-optimize" approach in order to compare it with our results.}

In conclusion, as another goal it might be interesting to study the performance of new stabilization techniques for the online phases, such as the Online Vanishing Viscosity and the Online Rectification methods \cite{ali2018stabilized, chakir2019non, maday1989analysis}. {Moreover, the extension of this setting to the uncertainty certification context will be the topic of future research.}

\section*{Acknowledgements}
We acknowledge the support by European Union Funding for Research and Innovation -- Horizon 2020 Program -- in the framework of European Research Council Executive Agency: Consolidator Grant H2020 ERC CoG 2015 AROMA-CFD project 681447 ``Advanced Reduced Order Methods with Applications in Computational Fluid Dynamics''. We also acknowledge the PRIN 2017  ``Numerical Analysis for Full and Reduced Order Methods for the efficient and accurate solution of complex systems governed by Partial Differential Equations'' (NA-FROM-PDEs) and the INDAM-GNCS project ``Tecniche Numeriche Avanzate per Applicazioni Industriali''.
The computations in this work have been performed with RBniCS \cite{RBniCS} library, developed at SISSA mathLab, which is an implementation in FEniCS \cite{logg2012automated} of several reduced order modelling techniques; we acknowledge developers and contributors to both libraries.

\bibliographystyle{plain}
\bibliography{BIB.bib}

\end{document}